%% file: conicbinary.tex
\documentclass[opre,copyedit]{informs3}

\OneAndAHalfSpacedXI %

\usepackage{natbib}
\bibpunct[, ]{(}{)}{,}{a}{}{,}%
\def\bibfont{\small}%

\usepackage{tikz,subfigure,graphicx,transparent,caption,color,multirow,etoolbox}

\newcommand{\cF}{\mathcal{F}}

\usepackage{hyperref}
\usepackage{mathtools}
\usepackage{bm}

\newcommand{\fullrank}[1]{}
\definecolor{gold}{rgb}{0.85,0.65,0}

\def\iff{\text{if and only if}}

\renewcommand{\theequation}{\mbox{\arabic{equation}}}

\def\fnote#1{\footnote}
\newcommand{\epr}{\hfill\hbox{\hskip 4pt \vrule width 5pt height 6pt depth 1.5pt}\vspace{0.0cm}\par}

\newcommand{\norm}[1]{\ensuremath{\left\lVert #1 \right\rVert}}

\def\C{{\mathbb{C}}}

\def\K{{\mathbb{K}}}
\def\L{{\mathbb{L}}}

\def\P{{\mathbb{P}}}

\def\R{{\mathbb{R}}}
\def\S{{\mathbb{S}}}

\def\Z{{\mathbb{Z}}}
\def\Le{{\mathbb{L}}}

\def\cF{{\cal F}}
\def\cG{{\cal G}}
\def\cH{{\cal H}}

\def\cM{{\cal M}}

\def\cO{{\cal O}}
\def\cP{{\cal P}}
\def\cQ{{\cal Q}}
\def\cR{{\cal R}}
\def\cS{{\cal S}}

\def\cW{{\cal W}}
\def\cX{{\cal X}}

\DeclareMathOperator{\Proj}{Proj}

\DeclareMathOperator{\Diag}{Diag}

\DeclareMathOperator{\Epi}{epi}
\DeclareMathOperator{\epi}{epi}

\DeclareMathOperator{\conv}{conv}

\DeclareMathOperator{\cone}{cone}

\DeclareMathOperator{\rec}{rec}

\def\log{\mathop{{\rm log}}}

\graphicspath{{fig/}}

\TheoremsNumberedThrough     %
\ECRepeatTheorems

\EquationsNumberedThrough    %

\begin{document}

	\RUNAUTHOR{K{\i}l{\i}n\c{c}-Karzan, K\"{u}\c{c}\"{u}kyavuz, Lee,  Shafieezadeh-Abadeh}
	
	\RUNTITLE{Conic Mixed-Binary Sets} %
	
	\TITLE{Conic Mixed-Binary Sets: Convex Hull Characterizations and Applications}
	
	\ARTICLEAUTHORS{%
		\AUTHOR{Fatma K{\i}l{\i}n\c{c}-Karzan}
		\AFF{Tepper School of Business, Carnegie Mellon University, Pittsburgh, PA 15213, USA, \EMAIL{fkilinc@andrew.cmu.edu}} %
		\AUTHOR{Simge K\"{u}\c{c}\"{u}kyavuz}
		\AFF{Department of Industrial Engineering and Management Sciences, Northwestern University, Evanston, IL 60208, USA, \EMAIL{simge@northwestern.edu}}
		\AUTHOR{Dabeen Lee}
		\AFF{Department of Industrial and Systems Engineering, KAIST, Daejeon 34141, South Korea, \EMAIL{dabeenl@kaist.ac.kr}}
		\AUTHOR{Soroosh Shafieezadeh-Abadeh}
		\AFF{Tepper School of Business, Carnegie Mellon University, Pittsburgh, PA 15213, USA, \EMAIL{sshafiee@andrew.cmu.edu}} %
	} %

	\input{abstract.tex}
	
	\maketitle

	\input{introduction}

	\input{conic_set}
	\input{extension}
	\input{application}

	\ACKNOWLEDGMENT{We thank the review team for their constructive feedback that lead to improved presentation of the content. This research is supported, in part, by ONR grants N00014-19-1-2321 and N00014-22-1-2602, AFOSR grant FA9550-22-1-0365, the Institute for Basic Science (IBS-R029-C1, Y2), the FOUR Brain Korea 21 Program (NRF-5199990113928), the National Research Foundation of Korea (NRF-2022M3J6A1063021), the KAIST Starting Fund (KAIST-G04220016), NSF grant CMMI 1454548, and Early Postdoc Mobility Fellowship SNSF grant P2ELP2\_195149.}

\input{appendix}
	
	\bibliographystyle{plainnat} 
	\bibliography{mybibfile}

\end{document}

%% file: abstract.tex
\ABSTRACT{
We consider a general conic mixed-binary set where each homogeneous conic constraint $j$ involves an affine function of
independent continuous variables and an epigraph variable associated with a nonnegative function, $f_j$, of common
binary variables. Sets of this form naturally arise as substructures in a number of applications including mean-risk
optimization, chance-constrained problems, portfolio optimization, lot-sizing and scheduling, fractional programming, variants of the best subset selection problem, a class of sparse semidefinite programs, 
and distributionally robust chance-constrained programs. 
%When all of the functions $f_j$'s are submodular, we give a convex hull description of this set that relies on characterizing the epigraphs of $f_j$'s. 
We give a convex hull description of this set that relies on simultaneous characterization of the epigraphs of $f_j$'s, which is easy to do when all functions $f_j$'s are submodular.
Our result unifies and  generalizes an existing result in two important directions. First, it considers \emph{multiple
	general convex cone} constraints instead of a single second-order cone type constraint. Second, it takes \emph{arbitrary nonnegative functions} instead of  a specific submodular function obtained from the square root of an affine function. 
We close by demonstrating the applicability of our results in the context of a number of %
problem classes.
}

\KEYWORDS{Conic mixed-binary sets, conic quadratic optimization, convex hull, submodularity, fractional binary optimization, best subset selection}

%% file: introduction.tex
\section{Introduction}\label{sec:intro}

In this paper we study the following set 
\begin{align}\label{eq:conicSet}
\cS(f,\K)\coloneqq\left\{(x,z)\in\R^m\times\{0,1\}^n:~~ \exists y\in\R \text{ s.t. } y\geq f(z),~Ax +By\in \K   \right\},
\end{align}
where $f:\{0,1\}^n\to\R_+$ is a nonnegative function, $\K$ is a closed convex pointed cone, %
and $A,B$ are any matrices of appropriate dimension. 
Throughout, for $a\in \Z_+$, we let $[a]\coloneqq\{1,\dots,a\}$.

The set $\cS(f,\K)$ we consider is a generalization of the two fundamental sets studied by \cite{atamturk2019conic-quad} defined as follows:
\begin{align}
\cH&\coloneqq \left\{(x,z)\in\R^{m}_+\times\{0,1\}^n:~\sqrt{\sigma+\sum_{i\in[n]}c_iz_i+ \sum_{j\in[m-1]} d_j x_j^2 }\leq x_{m} \right\},\label{H}\\
\cR&\coloneqq \left\{(x,z)\in\R^{m}_+\times\{0,1\}^n:~\sigma+\sum_{i\in[n]}c_iz_i+ \sum_{j\in[m-2]} d_j x_j^2 \leq 4x_{m-1}x_{m}   \right\},\label{R}
\end{align}
where $c\in\R_+^n$, $d\in\R_+^{m-1}$, and $\sigma\in\R_+$. \cite{atamturk2019conic-quad} demonstrate that the set $\cH$ naturally arises in mean-risk minimization and chance-constrained programs, and the set $\cR$ appears as a substructure in robust conic quadratic interdiction, lot-sizing and scheduling, queueing system design, binary linear fractional problems, Sharpe ratio maximization, and portfolio optimization. The main contribution of \cite{atamturk2019conic-quad} is deriving the convex hulls of $\cH$ and $\cR$.

In Section~\ref{sec:convex_hull}, we prove that the convex hull of the set $\cS(f,\K)$ for any nonnegative function $f$ and arbitrary  closed convex pointed cone $\K$ \fullrank{and a full column rank linear map $A$} is given by %
\begin{equation}\label{eq:shat-def}
\widehat{\cS}(f,\K)\coloneqq\left\{(x,z)\in\R^m\times[0,1]^n:~~\exists y\in\R_+ \text{ s.t. } (y,z)\in\conv(\Epi(f)),~ Ax +By\in \K   \right\},
\end{equation}
where  given a function $f$,  $\Epi(f)$ denotes its epigraph, i.e., $\Epi(f)\coloneqq \{(y,z)\in\R\times\{0,1\}^n: y\geq f(z) \}$ and $\conv(\Epi(f))$ denotes its convex hull. Therefore, the convex hull of $\cS(f,\K)$ is given by the inequalities describing $\conv(\Epi(f))$ and the homogeneous conic constraint $Ax +By \in \K$. This characterization  highlights, in particular, that the challenge in developing the convex hull of $\cS(f,\K)$ is solely determined by the complexity of the convex hull characterization of the epigraph of the function $f$. 
 In fact, we prove a more general result that covers this characterization as a special case. In particular, we consider a set with \emph{multiple} conic constraints, each of which involves an affine function of
independent continuous variables and an epigraph variable associated with a nonnegative function of \emph{common}
binary variables.

 There are a number of cases where characterizing %
 $\conv(\Epi(f))$ is easy. We discuss a few such cases at the end of Section~\ref{sec:convex_hull}.  Most notably, %
 when $f$ is a nonnegative submodular function, imposing the condition that $(y,z)\in\conv(\Epi(f))$ is equivalent to applying the \emph{extended polymatroid inequalities}~\cite[]{atamturk2008polymatroids,lovasz1983submodular}.   For completeness, in Appendix~\ref{sec:inequalities} we %
review these inequalities. % and their generalization to arbitrary set functions, namely the polar inequalities, \citep{atamturk2020polar}. %
In Section~\ref{sec:preliminaries}, we also discuss how our results can be applied to more general functions and non-homogeneous conic constraints. 

The set $\cS(f,\K)$ naturally arises in a number of other applications including fractional binary optimization, best subset selection, a class of sparse semidefinite programs, and distributionally robust optimization. We discuss these in Section~\ref{sec:applications}. %and Appendix~\ref{sec:DRO-CCP}. 
For example, %
the sets $\cH$ and $\cR$ studied by \cite{atamturk2019conic-quad} are indeed special cases of $\cS(f,\K)$ where $\K$ is taken to be the direct product of a second-order cone (SOC) and the nonnegative orthant, and the function $f(z)$ is restricted to be the square root of an affine function of $z$. We elaborate  this connection in Section~\ref{sec:soc}. Consequently, even in the case of a single conic constraint and a single function $f$, our convex hull characterization immediately generalizes the results from \cite{atamturk2019conic-quad} in two directions%
: (1) by allowing general closed convex pointed cones $\K$ as opposed to the standard SOC, and (2) by allowing a general nonnegative function $f(z)$ as opposed to the specific one studied by \cite{atamturk2019conic-quad}. We discuss applications of our framework with \emph{multiple} conic constraints and \emph{multiple} functions on \emph{common} binary variables in Section~\ref{sec:soc-multiple}. In particular, this generalization allows us to cover the fractional optimization models that appear in a broad range of applications including %
multinomial logit (MNL) choice models in assortment optimization, set covering, market share based facility location, stochastic service systems, bi-clustering, and optimization of boolean query for databases.
We discuss application of our framework to the best subset selection problem in machine learning in Section~\ref{sec:BSS}. In this application, the function $f$ is an exponential function, which demonstrates the applicability of our result beyond the square root function.  
Furthermore, our results in this context pave the way to use standard solvers for this problem all the while exploiting the submodular structure as opposed to the approach of~\cite{gomez2020fractional} based on parameterizing the fractional objective function and applying a customized Newton-type method. In Appendix~\ref{app:numerical} we provide a numerical study illustrating the effectiveness of our general purpose approach against this custom made method for this specific problem from Section~\ref{sec:BSS}. In Section~\ref{sec:misdp} we consider sparse approximation of the sum of positive semidefinite matrices as an application, which fits to our framework by taking $\K$ to be the positive semidefinite cone. 
 We also consider an application in distributionally robust chance-constrained programming under Wasserstein ambiguity in Appendix~\ref{sec:DRO-CCP} and show that our results can be used to strengthen the mixed-integer conic reformulation for the case of mixed-binary decision variables, thereby generalizing an existing strengthening that assumes pure binary decisions in the original chance constraint.

%% file: conic_set.tex
\section{Convex hull characterization}\label{sec:convex_hull}
  
We next examine a generalization of our set defined by
\begin{align}\label{eq:S-multiple}
\cS(\{f_j\}_{j\in[p]},\{\K_j\}_{j\in[p]}) &\coloneqq  \left\{(x,z)\in\R^{mp}\times\{0,1\}^{n}:~\exists y\in\R^p \text{ s.t. } \begin{array}{l} y_j\geq f_j(z),~~\forall j\in[p], \\ A^jx^j + B^j y_j \in \K_j,~~\forall j\in[p] \end{array} \right\},
\end{align}
where the functions $f_j:\{0,1\}^n\to\R_+$ for $j\in[p]$ are nonnegative functions and $\K_j$ for $j\in[p]$ are closed  convex pointed cones. $A^j$ and $B^j$ for $j\in[p]$ are matrices of appropriate dimensions. The lengths of the continuous vectors $x^1,\ldots,x^p$ may be different, but we focus on the setting of equal lengths for simplicity. Our theoretical developments can be easily extended to the case of heterogeneous lengths. 

Note that the following is a convex relaxation of $\cS(\{f_j\}_{j\in[p]},\{\K_j\}_{j\in[p]})$:
\begin{eqnarray}
&\widehat{\cS}(\{f_j\}_{j\in[p]},\{\K_j\}_{j\in[p]}) \coloneqq  
\left\{(x,z)\in\R^{mp}\times[0,1]^{n}:~\exists y\in\R^p \text{ s.t. }~\begin{array}{l} (y,z)\in\conv\left(\cG\left(\{f_j\}_{j\in[p]}\right)\right),\\  A^jx^j + B^j y_j \in \K_j,~\forall j\in[p] \end{array}  \right\}, \\ \label{eq:Shat-multiple}
&  
\text{where}\qquad\qquad\cG\left(\{f_j\}_{j\in[p]}\right) \coloneqq \left\{(y,z)\in\R^{p}\times\{0,1\}^{n}:~y_j\geq f_j(z),~~\forall j\in[p]\right\} . \label{eq:Gmultiple}
\end{eqnarray}
We note that the functions $f_j$ for $j\in[p]$ take the \emph{same} binary variables $z$. The constraint $y_j\geq f_j(z)$ gives rise to the epigraph of $f_j$ for each $j$, and therefore, $\cG\left(\{f_j\}_{j\in[p]}\right)$ can be viewed as the ``intersection" of the epigraphs. Our main result establishes that, under minor assumptions, the convex hull of the set $\cS(\{f_j\}_{j\in[p]},\{\K_j\}_{j\in[p]})$ is indeed given precisely by $\widehat{\cS}(\{f_j\}_{j\in[p]},\{\K_j\}_{j\in[p]})$.

To prove our convex hull result, we first establish a technical proposition that applies to more general mixed-integer sets. We consider sets of the form
\begin{align}\label{eq:cQ}
\cQ(\cG,\{\K_j\}_{j\in[p]})&\coloneqq  \left\{(x,y,z)\in\R^{mp}\times\R^p\times\{0,1\}^{n}:~~(y,z)\in\cG,~~ A^jx^j + B^j y_j \in \K_j,~\forall j\in[p]  \right\},
\end{align}
where $\cG\subseteq \R^{p}_+\times\{0,1\}^{n}$ and $\K_j$'s for $j\in[p]$ are closed  convex pointed cones. We further assume that $\conv(\cG)$ is a closed pointed set such that %
$\cG|_{z=\bar z}\coloneqq\cG\cap \{(y,z): z=\bar z\}$ for any fixed $\bar z\in\{0,1\}^n$ is convex, which means that $\cG|_{z=\bar z}$ is the face of $\conv(\cG)$ defined by $z=\bar z$.
Note that $\cG\left(\{f_j\}_{j\in[p]}\right)$ is indeed contained in $\R^{p}_+\times\{0,1\}^{n}$, as $f_j$ for $j\in[p]$ are nonnegative functions. Moreover, $\cG\left(\{f_j\}_{j\in[p]}\right)|_{z=\bar z}$ for a fixed $\bar z\in\{0,1\}^n$ is defined by linear constraints only, so it is convex. Observe that $\cS(\{f_j\}_{j\in[p]},\{\K_j\}_{j\in[p]})$ is precisely the projection of $\cQ(\cG\left(\{f_j\}_{j\in[p]}\right),\{\K_j\}_{j\in[p]})$  onto the $(x,z)$-space.  
Given this relation, to derive our promised convex hull result for the set $\cS(\{f_j\}_{j\in[p]},\{\K_j\}_{j\in[p]})$ in the original space, we first establish the following result for sets of the form $\cQ(\cG,\{\K_j\}_{j\in[p]})$ defined in the extended space.

\begin{theorem}\label{thm:conv-trick}
Consider $\cQ(\cG,\{\K_j\}_{j\in[p]})$ defined as in~\eqref{eq:cQ} for some set $\cG\subseteq \R^{p}_+\times\{0,1\}^{n}$. Then, $\conv(\cQ(\cG,\{\K_j\}_{j\in[p]})) \subseteq  \widehat{\cQ}(\cG,\{\K_j\}_{j\in[p]})$, where 
\begin{align}\label{eq:cQhat}
\widehat{\cQ}(\cG,\{\K_j\}_{j\in[p]}) \coloneqq \!\! \left\{(x,y,z)\in\R^{mp}\times\R^p\times[0,1]^{n}:(y,z)\in\conv(\cG),~A^jx^j + B^j y_j \in \K_j,\,\forall j\in[p]  \right\}.\!\!
\end{align}
Moreover, if $\K_j$ is a closed  convex pointed cone for each $j\in[p]$ and $\cG$ is such that $\conv(\cG)$ is closed and $\cG|_{z=\bar z}$ for any fixed $\bar z\in\{0,1\}^n$ is convex,
then 
$
\conv(\cQ(\cG,\{\K_j\}_{j\in[p]})) =  \widehat{\cQ}(\cG,\{\K_j\}_{j\in[p]}).
$
\end{theorem}
\proof{{\bf Proof.%
}}
For brevity, define $\cQ\coloneqq \cQ(\cG,\{\K_j\}_{j\in[p]})$ and $\widehat{\cQ}\coloneqq \widehat{\cQ}(\cG,\{\K_j\}_{j\in[p]})$. 
Note that $\conv(\cQ) \subseteq \widehat{\cQ}$ is immediate because $\cQ \subseteq \widehat{\cQ}$ and $\widehat{\cQ}$ is convex, so it is sufficient to show that $\widehat{\cQ}\subseteq\conv(\cQ)$ holds under the stated conditions. 
	
First, note that $\widehat{\cQ}$ is closed since it is the intersection of closed convex sets defined by $\conv(\cG)$ (which is closed by the premise of the theorem) and each of the constraints  $A^jx^j+B^jy_j\in\K_j$ defines a closed set (as each cone $\K_j$ is closed). 
Let $L\coloneqq \!\! \left\{(x,0,0)\in\R^{mp}\times\R^p\times\R^{n}:  A^jx^j=0,\forall j\in[p]\right\}$. Then, $L$ is a linear subspace and 
$\cQ =\cQ +L$ (it is clear that $\cQ\subseteq \cQ+L$ and for the other direction we note that for any $(\bar x,\bar y,\bar z)\in\cQ$ and any $\tilde x \in L$, we have $(\bar y,\bar z)\in\cG$ and $A^j (\bar x + \tilde x)^j + B^j \bar{y}_j = A^j \bar{x}^j + B^j \bar{y}_j \in\K_j$ for all $j\in[p]$ implying  $(\bar x+\tilde{x},\bar y,\bar z)\in\cQ$.  Similarly, we also see that $\widehat{\cQ}=\widehat{\cQ}+L$.

We define $\widetilde{\cQ}:=\widehat{\cQ}\cap L^\perp$ where $L^\perp$ is the linear subspace orthogonal to $L$. Then, $\widehat{\cQ} = \widetilde{\cQ} + L$. % {\crd (maybe add more details)}.  
Note that $\widetilde{\cQ}$ is closed and convex because it is the intersection of two closed convex sets, namely $\widehat{\cQ}$ and the linear subspace $L^\perp$. We further claim that $\widetilde{\cQ}$ is pointed. Assume for contradiction that $\widetilde{\cQ}$ is not pointed. Then, there exists $(d_x,d_y,d_z)\in(\R^{mp}\times\R^p\times\R^n)\setminus \{(0,0,0)\}$ such that $(x,y,z)+\alpha(d_x,d_y,d_z)\in \widetilde{\cQ}$ for any $(x,y,z)\in\widetilde{\cQ}$ and any $\alpha\in\R$. As $\cG\subseteq \R_+^p\times \{0,1\}^n$, $\conv(\cG)$ is pointed and thus from $\widetilde{\cQ} \subseteq (\R^{mp}\times\conv(\cG))$ we deduce that $d_y=0$ and $d_z=0$. Then, we must have $d_x\neq0$ as $(d_x,d_y,d_z)$ is nonzero. Moreover, as $\widetilde{\cQ} \subseteq \widehat{\cQ}$, we deduce that for all $\alpha\in\R$, $(x+\alpha d_x,y,z) \in \widehat{\cQ}$ and thus $A^j (x+\alpha d_x)^j + B^j y_j \in\K_j$ for all $j\in[p]$. Since each $\K_j$ is a pointed cone, this is possible if and only if $A^j d_x^j =0$ for all $j\in[p]$, i.e., $(d_x,0,0)\in L$. On the other hand, by definition of $\widetilde{\cQ}$, we have $\widetilde{\cQ} \subseteq L^\perp$ and thus $(x+\alpha d_x,y,z) \in L^\perp$ for all $\alpha \in\R$, which is possible if and only if $(d_x,0,0)\in L^\perp$. So, we have $(d_x,0,0)\in L \cap L^\perp=\{(0,0,0)\}$ which contradicts $d_x \neq 0$. Therefore, we conclude that $\widetilde{\cQ}=\widehat{\cQ}\cap L^\perp$ is a closed convex and pointed set.

As $\widetilde{\cQ}$ is a closed convex pointed set, by defining $E$ and $R$ as its sets of extreme points and extreme rays, respectively, we can write $\widetilde{\cQ}=\conv(E)+\cone(R)$ and thus $\widehat{\cQ} = \widetilde{\cQ} + L=\conv(E)+\cone(R)+L$. To prove that $\widehat{\cQ}\subseteq \conv(\cQ)$, it suffices to show that $E\subseteq \cQ$, $\cone(R) \subseteq \rec(\conv(\cQ))$ and that $L\subseteq \rec(\conv(\cQ))$, where $\rec(\conv(\cQ))$ denotes the recessive cone of $\conv(\cQ)$. Recall, $\cQ=\cQ +L$, and thus $\conv(\cQ)=\conv(\cQ+L)=\conv(\cQ)+L$. Hence, $L\subseteq\rec(\conv(\cQ))$.

Consider any $(d_x, d_y, d_z)\in \cone(R)$. Then, $A^j d_x^j + B^j (d_y)_j\in\K_j$ for $j\in[p]$ and  $(d_y, d_z)\in\rec(\conv(\cG))$. Assume for contradiction that $(d_x, d_y, d_z)\not\in\rec(\conv(\cQ))$. Then, there exists $(\bar x, \bar y, \bar z)\in\conv(\cQ)$ such that $(\bar x, \bar y, \bar z)+(d_x, d_y, d_z)\not\in \conv(\cQ)$. As $(\bar x, \bar y, \bar z)\in\conv(\cQ)$, it can be written as a convex combination of some points, denoted $(x^i,y^i,z^i)$ for $i\in I$, in $\cQ$. Hence, $(\bar x, \bar y, \bar z)=\sum_{i\in I}\alpha_i(x^i,y^i,z^i)$ for some $\alpha\geq0$ with $\sum_{i\in I}\alpha_i=1$. Since $(x^i,y^i,z^i)\in \cQ$, we have $A^j (x^i)^j + B^j (y^i)_j\in\K_j$ for $j\in[k]$, and therefore, $A^j (x^i+d_x)^j + B^j (y^i+d_y)_j\in\K_j$ for $j\in[k]$. Since $(d_y, d_z)\in\rec(\conv(\cG))$, $(y^i+d_y, z^i+d_z)\in \conv(\cG)$. Moreover, $\conv(\cG)\subseteq \R_+^p\times [0,1]^n$ implies that $d_z=0$, so $(y^i+d_y, z^i+d_z)\in \conv(\cG)\cap\left\{(y,z):z=z^i\right\}=\cG|_{z=z^i}\subseteq \cG$. Therefore, $(x^i+d_x,y^i+d_y,z^i+d_z)\in \cQ$ for $i\in I$, which implies that $(\bar x, \bar y, \bar z)+(d_x, d_y, d_z)=\sum_{i\in I}\alpha_i(x^i+d_x,y^i+d_y,z^i+d_z)\in \conv(\cQ)$, a contradiction. Therefore, $(d_x, d_y, d_z)\in\rec(\conv(\cQ))$, so $\cone(R) \subseteq \rec(\conv(\cQ))$.

Let $(\hat x, \hat y, \hat z)\in E$. Then, $(\hat{y},\hat{z})\in\conv(\cG)$, $A^j \hat{x}^j+B^j \hat{y}_j\in \K_j$  for all $j\in[p]$ and from $E\subseteq \widetilde{\cQ} \subseteq L^\perp$ we deduce $(\hat x, \hat y, \hat z)\in L^\perp$. Also, as $\conv(\cG)\subseteq\R^p_+\times[0,1]^n$, $\hat y\in\R^p_+$. We claim that $(\hat{y},\hat{z})$ must be in $\cG$. We will prove this by showing that $(\hat{y},\hat{z})$ must be an extreme point of $\conv(\cG)$. Assume for contradiction that there exist distinct points $(\bar{y},\bar{z}) \in \conv(\cG)$ and $(\tilde{y},\tilde{z})\in\conv(\cG)$ such that $(\hat{y},\hat{z}) = {1\over 2} (\bar{y},\bar{z}) + {1\over 2} (\tilde{y},\tilde{z})$. From $\conv(\cG)\subseteq\R^p_+\times[0,1]^n$, we deduce that $\bar{y},\tilde{y}\in\R^p_+$. Moreover, if for some index $j\in[p]$ we have $\hat{y}_j=0$, we must also have $\bar{y}_j=\tilde{y}_j=0$. For each $j\in[p]$, define $\bar{x}^j \coloneqq { \bar{y}_j \over \hat{y}_j} \hat{x}^j$ whenever $\hat{y}_j > 0$, and $\bar{x}^j \coloneqq \hat{x}^j$ whenever $\hat{y}_j = 0$.  Note that for all $j\in[p]$, by definition $\bar{x}^j$ is a positive multiple of $\hat{x}^j$ and as $(\hat{x},\hat{y},\hat{z})\in L^\perp$ from the definition of $L$ we deduce that $(\bar{x},\bar{y},\bar{z})\in L^\perp$ as well. Consider any $j\in[p]$. Then, when $\hat{y}_j = 0$, we have $\bar{y}_j=0$ as well as $\bar{x}^j=\hat{x}^j$, and thus $(\bar{x}^j,\bar{y}_j)=(\hat{x}^j,\hat{y}_j)$ and hence $A^j \bar{x}^j+B^j \bar{y}_j=A^j \hat{x}^j+B^j \hat{y}_j\in \K_j$. When $\hat{y}_j > 0$, by definition we have $(\bar{x}^j,\bar{y}_j)= { \bar{y}_j \over \hat{y}_j} (\hat{x}^j,\hat{y}_j)$. Moreover, because  $A^j \hat{x}^j+B^j \hat{y}_j\in \K_j$ holds, ${ \bar{y}_j \over \hat{y}_j}\in\R_+$ and $\K_j$ is a closed convex pointed cone, we deduce that $A^j \bar{x}^j+B^j \bar{y}_j\in \K_j$ as well. Hence, in either case we conclude that $(\bar{x},\bar{y},\bar{z}) \in \widetilde{\cQ}$. Similarly, for each $j\in[p]$, define $\tilde{x}^j \coloneqq { \tilde{y}_j \over \hat{y}_j} \hat{x}^j$ whenever $\hat{y}_j > 0$, and $\tilde{x}^j \coloneqq \hat{x}^j$ whenever $\hat{y}_j = 0$. As before, we deduce that $(\tilde{x},\tilde{y},\tilde{z}) \in \widetilde{\cQ}$. Finally, for each $j\in[p]$ such that $\hat{y}_j>0$ we deduce \[{1\over 2} (\bar{x}^j + \tilde{x}^j )= {1\over 2} \left({\bar{y}_j \over \hat{y}_j} \hat{x}^j + {\tilde{y}_j \over \hat{y}_j} \hat{x}^j \right) = {\bar{y}_j + \tilde{y}_j \over 2 \hat{y}_j}  \hat{x}^j = \hat{x}^j,\] where the last equation follows from $(\hat{y},\hat{z}) = {1\over 2} (\bar{y},\bar{z}) + {1\over 2} (\tilde{y},\tilde{z})$. Also, for each $j\in[p]$ such that $\hat{y}_j=0$, we have $\bar{x}^j = \tilde{x}^j=\hat{x}^j$. Thus, $(\hat x, \hat y, \hat z)= {1\over 2} (\bar{x},\bar{y},\bar{z}) + {1\over 2} (\tilde{x},\tilde{y},\tilde{z})$. This contradicts $(\hat x, \hat y, \hat z)\in E$. Therefore,  $(\hat{y},\hat{z})$ is an extreme point of $\conv(\cG)$, and thus we must have $(\hat{y},\hat{z})\in\cG$. Hence, we have shown that $E\subseteq \cQ$, as required.
\Halmos
\endproof

Theorem~\ref{thm:conv-trick} is instrumental in proving the following main theorem that gives the convex hull characterization of $\cS(\{f_j\}_{j\in[p]},\{\K_j\}_{j\in[p]})$.
 
\begin{theorem}\label{thm:multipleFunctions}
For each $j\in[p]$, let $f_j:\{0,1\}^n\to\R_+$ be a nonnegative %
function and $\K_j$ be a closed  convex pointed cone. %
Then, \fullrank{$\conv\left(\cS(\{f_j\}_{j\in[p]},\{\K_j\}_{j\in[p]})\right) \subseteq \widehat{\cS}(\{f_j\}_{j\in[p]},\{\K_j\}_{j\in[p]})$ %
defined as in \eqref{eq:S-multiple} and \eqref{eq:Shat-multiple}. Moreover, when the matrix $A^j$ %
 	has full column rank for each $j\in[p]$, we have} $\conv\left(\cS(\{f_j\}_{j\in[p]},\{\K_j\}_{j\in[p]})\right) = \widehat{\cS}(\{f_j\}_{j\in[p]},\{\K_j\}_{j\in[p]})$.
\end{theorem}
\proof{{\bf Proof.%
}}
First, we observe that $\cG\coloneqq \cG\left(\{f_j\}_{j\in[p]}\right)\subseteq\R^p_+\times\{0,1\}^n$ as $f_j$'s are nonnegative functions and that 
$ %
\cS(\{f_j\}_{j\in[p]},\{\K_j\}_{j\in[p]}) = \Proj_{x,z}\left(\cQ (\cG,\{\K_j\}_{j\in[p]})\right).
$ %
Moreover, $\conv(\cG)$ is a polyhedral set, and thus $\cG|_{z=\bar z}$ for any $z=\bar z$ is convex. 
Next, %
recall that the convex hull and projection operations commute, i.e., the projection of the convex hull of a set is equal to the convex hull of the projection of a set \cite[see e.g.,][Exercise 3.34]{CCZ2014book}. Therefore,
\begin{align*}
\conv\left(\cS(\{f_j\}_{j\in[p]},\{\K_j\}_{j\in[p]}) \right) &= \conv\left(\Proj_{x,z}\left(\cQ (\cG,\{\K_j\}_{j\in[p]})\right) \right) \\
& = \Proj_{x,z}\left(\conv\left(\cQ (\cG,\{\K_j\}_{j\in[p]})\right) \right)\\
&\subseteq \Proj_{x,z} \left(\widehat{\cQ} (\cG,\{\K_j\}_{j\in[p]})\right)  
 = \widehat{\cS}(\{f_j\}_{j\in[p]},\{\K_j\}_{j\in[p]}),   
\end{align*}
where the last equality holds from 
$ %
\widehat{\cS}(\{f_j\}_{j\in[p]},\{\K_j\}_{j\in[p]}) = \Proj_{x,z}\left(\widehat{\cQ} (\cG,\{\K_j\}_{j\in[p]})\right).
$ %

Because $\cG\subseteq\R^p_+\times\{0,1\}^n$, $\conv(\cG)$ %
has at least one extreme point as well. This 
\fullrank{together with the premise that  the matrix $A^j$  is of full column rank for each $j\in[p]$}  
then satisfies all conditions of the last part of 
Theorem~\ref{thm:conv-trick} which implies $\conv\left(\cQ (\cG,\{\K_j\}_{j\in[p]})\right)=\widehat{\cQ} (\cG,\{\K_j\}_{j\in[p]})$. Thus, we deduce that $\conv\left(\cS(\{f_j\}_{j\in[p]},\{\K_j\}_{j\in[p]})\right) = \widehat{\cS}(\{f_j\}_{j\in[p]},\{\K_j\}_{j\in[p]})$ holds.
% under the stated conditions. %
\Halmos
\endproof

The interesting applications of Theorem~\ref{thm:multipleFunctions} arise whenever it is relatively easy to give the explicit convex hull description of the intersections of the epigraphs of nonnegative functions $f_j:\{0,1\}^n\to\R_+$. One such case is when $f_j$'s are nonnegative and submodular, where $\conv(\Epi(f_j))$'s are described by the \emph{extended polymatroid inequalities}. When $f_j$'s are general set functions, a recent work of~\cite{atamturk2020polar} discusses the \emph{polar inequalities} that generalize the extended polymatroid inequalities. Hence, the polar inequalities are valid for $\cS(\{f_j\}_{j\in[p]},\{\K_j\}_{j\in[p]})$ for any set of nonnegative functions $f_j$'s. %
We briefly review these inequalities %
in Appendix~\ref{sec:inequalities}.

Next, we focus on the $p=1$ case, introduced in Section~\ref{sec:intro}, where we consider a single set function $f$ and a  closed  convex pointed cone $\K$. To simplify our notation in this case, we  define 
 \[
\cS(f,\K)\coloneqq  \cS\left(\{f\},\{\K\}\right)=\left\{(x,z)\in\R^{m}\times\{0,1\}^{n}:~\exists y\in\R\text{ s.t. } y\geq f(z),~Ax + B y \in \K\right\}.
\] 
 Consider the set
\begin{equation}\label{S-single}
\overline\cS(f,\K)\coloneqq  \left\{(x,z)\in\R^{m}\times\{0,1\}^{n}:~Ax + B f(z) \in \K\right\}.
\end{equation}
Then, clearly $\overline\cS(f,\K)\subseteq\cS(f,\K)$, and we arrive at the following immediate corollary of Theorem~\ref{thm:multipleFunctions}.

\begin{corollary}\label{cor:singleFunction}
Let $f:\{0,1\}^n\to\R_+$ be a nonnegative function and $\K$ be a closed  convex pointed cone %
Then, $\conv(\overline\cS(f,\K))\subseteq \widehat{\cS}(f,\K)$. If $\overline{\cS}(f,\K)=\cS(f,\K)$ further holds\fullrank{and the matrix $A$ has full column rank}, then 
$ %
\conv(\overline\cS(f,\K))=\widehat{\cS}(f,\K).
$ %
\end{corollary}
\proof{{\bf Proof.%
}}
Since $\overline\cS(f,\K)\subseteq \cS(f,\K)$, it follows that $\conv(\overline\cS(f,\K)) \subseteq \conv(\cS(f,\K))\subseteq \widehat{\cS}(f,\K)$,
where the last containment follows from  
Theorem~\ref{thm:multipleFunctions} applied to the $p=1$ case. %
If we have $\overline{\cS}(f,\K)=\cS(f,\K)$, then $\conv(\overline\cS(f,\K))=\conv(\cS(f,\K))$. Moreover, \fullrank{when the matrix $A$ %
	is of full column rank,} %
Theorem~\ref{thm:multipleFunctions} implies that $\conv(\cS(f,\K))= \widehat{\cS}(f,\K)$, and thus $\conv(\overline\cS(f,\K))=\widehat{\cS}(f,\K)$ indeed holds. %
\Halmos
\endproof

In general, the condition $\overline{\cS}(f,\K)=\cS(f,\K)$ does not always hold. 
For example, suppose that $n=1$, $\K=\R_+$, and $f(z)=z$ which is %
	nonnegative over $\{0,1\}$. Then consider the corresponding set $\overline\cS(f,\K)=\left\{(x,z)\in\R\times\{0,1\}:\; -x+f(z)\geq0\right\}=\left\{(x,z)\in\R\times\{0,1\}: -x+z\geq0\right\}$. In this case, the associated set $\cS(f,\K)=\left\{(x,z)\in\R\times\{0,1\}:\; \exists y\in\R, -x+y\ge 0,\, y\geq z \right\}=\R\times\{0,1\}$. While $\overline\cS(f,\K)=\left\{(x,z)\in\R\times\{0,1\}: -x+z\geq0\right\}$, we have $\cS(f,\K)=\R\times\{0,1\}$. In particular, the point $(\bar x,\bar z)=(10,0)$ is in $\cS(f,\K)\setminus \overline \cS(f,\K)$ and thus the condition  $\overline{\cS}(f,\K)=\cS(f,\K)$ is not satisfied in this example.
This in turn indicates that $\overline\cS(f,\K)$ is not equal to $\widehat{\cS}(f,\K)$ in general. 

Nevertheless, there are some important examples, wherein the condition $\overline{\cS}(f,\K)=\cS(f,\K)$ indeed holds and thus $\conv(\overline\cS(f,\K))=\widehat{\cS}(f,\K)$. 
The following proposition provides a necessary and sufficient condition for $\conv(\overline\cS(f,\K))=\widehat{\cS}(f,\K)$.
\begin{proposition}\label{Condition*}
$\overline{\cS}(f,\K)=\cS(f,\K)$ if and only if the following condition is satisfied:
\begin{equation}\label{eq:Condition*}
\text{every $(x,z)$ satisfying $Ax + Bf(z)\in\K$ also satisfies $A(\alpha x) + Bf(z)\in \K$ for any $\alpha\geq 1$.}\tag{$\star$}
\end{equation}
\end{proposition}
\proof{{\bf Proof.%
}}
$\bm{(\Rightarrow)}$ Take $(x,z)$ such that $Ax + Bf(z)\in\K$. By definition, we have $(x,z)\in \overline\cS(f,\K)$. Then, by  $\overline{\cS}(f,\K)=\cS(f,\K)$, we have $(x,z)\in \cS(f,\K)$, i.e., there exists $y\in\R$ such that $Ax+By\in\K$ and $y\geq f(z)$. As $\K$ is a cone and $f$ is a nonnegative function, then for any $\alpha\ge1$, by letting $y'\coloneqq \alpha y$ we get $A(\alpha x)+B y' = A(\alpha x)+B(\alpha y)\in\K$ and $y'=\alpha y\geq y \geq f(z)$. Thus, by definition of $\cS(f,\K)$, we deduce $(\alpha x,z)\in \cS(f,\K)$ for any $\alpha\geq 1$. Now, since $(\alpha x,z)\in \cS(f,\K)$ for any $\alpha\geq 1$, once again using $\overline{\cS}(f,\K)=\cS(f,\K)$, we deduce $(\alpha x,z)\in \overline\cS(f,\K)$. Hence, $A(\alpha x) + Bf(z)\in \K$ for any $\alpha\geq 1$. 

$\bm{(\Leftarrow)}$ It suffices to prove that ${\cS}(f,\K)\subseteq\overline\cS(f,\K)$. To this end, take $(x,z)\in{\cS}(f,\K)$. Then there exists $y\geq0$ such that $y\geq f(z)$ and $Ax+By\in\K$. If $y=0$, then $f(z)=0$ as $f$ is a nonnegative function. This implies that $Ax+B f(z)=Ax+By\in\K$, so $(x,z)\in \overline\cS(f,\K)$. If $y>0$, it follows from $A(f(z)/y)x + Bf(z)=(f(z)/y)(Ax+By)\in \K$. In this case, since $y/f(z)\geq 1$, the premise of this direction implies that $Ax + Bf(z)=A(y/f(z))(f(z)/y)x + Bf(z)\in \K$, and therefore, $(x,z)\in\overline{\cS}(f,\K)$. Thus, ${\cS}(f,\K)\subseteq\overline\cS(f,\K)$, as required.
\Halmos
\endproof
We highlight a useful sufficient condition that implies Condition~\eqref{eq:Condition*}. %
\begin{remark}\label{rem:condition*}
Let $\K$ be a closed  convex pointed cone. %
Suppose $\K$, $A$ and $B$ are such that 
\begin{equation*}
\textup{for any $(x,z)$ satisfying $Ax + Bf(z)\in\K$, we also have $Ax\in\K$.}
\end{equation*}
For any $(x,z)$ with $Ax + Bf(z)\in\K$ and any $\alpha\geq 1$, we have $(\alpha-1)Ax\in\K$, so $A(\alpha x) + Bf(z) = (\alpha-1) Ax + \left(Ax + Bf(z) \right) \in \K$. Therefore, Condition~\eqref{eq:Condition*} is immediately satisfied.
\ifx\flagJournal\true \epr \fi
\end{remark}

We discuss in Section~\ref{sec:applications} some applications where Condition~\eqref{eq:Condition*} holds. In particular, using Remark~\ref{rem:condition*}, we will show in Section~\ref{sec:soc} that $\cH$ and $\cR$ defined in~\eqref{H} and~\eqref{R} satisfy Condition~\eqref{eq:Condition*} and consider other applications in Sections~\ref{sec:soc-multiple} -- \ref{sec:DRO-CCP}.

Corollary~\ref{cor:singleFunction} and Proposition~\ref{Condition*} can be extended to the case of multiple functions. Consider
\begin{equation}\label{eq:Sbar-multiple}
\overline\cS(\{f_j\}_{j\in[p]},\{\K_j\}_{j\in[p]}):=\left\{(x,z)\in\R^{mp}\times\{0,1\}^n:~A^jx^j+B^jf_j(z)\in \K_j,~~\forall j\in[p]\right\},
\end{equation}
which is a subset of $\cS(\{f_j\}_{j\in[p]},\{\K_j\}_{j\in[p]})$ defined in~\eqref{eq:S-multiple}.
\begin{corollary}\label{cor:multipleFunctions}
Let $f_j:\{0,1\}^n\to\R_+$ be a nonnegative function and $\K_j$ be a closed  convex pointed  cone %
for $j\in[p]$. Then, $\conv(\overline\cS(\{f_j\}_{j\in[p]},\{\K_j\}_{j\in[p]}))\subseteq \widehat{\cS}(\{f_j\}_{j\in[p]},\{\K_j\}_{j\in[p]})$. If $\overline{\cS}(\{f_j\}_{j\in[p]},\{\K_j\}_{j\in[p]})=\cS(\{f_j\}_{j\in[p]},\{\K_j\}_{j\in[p]})$ further holds\fullrank{and the matrix $A^j$ %
	has full column rank for each $j\in[p]$}, then
	$%
	\conv(\overline\cS(\{f_j\}_{j\in[p]},\{\K_j\}_{j\in[p]}))=\widehat{\cS}(\{f_j\}_{j\in[p]},\{\K_j\}_{j\in[p]}).
	$%
\end{corollary}
\proof{{\bf Proof.%
}}
Similar to the proof of Corollary~\ref{cor:singleFunction}.
\Halmos
\endproof
Moreover, we can characterize when $\overline{\cS}(\{f_j\}_{j\in[p]},\{\K_j\}_{j\in[p]})=\cS(\{f_j\}_{j\in[p]},\{\K_j\}_{j\in[p]})$.
\begin{proposition}\label{Condition**}
$\overline{\cS}(\{f_j\}_{j\in[p]},\{\K_j\}_{j\in[p]})=\cS(\{f_j\}_{j\in[p]},\{\K_j\}_{j\in[p]})$ if and only if the following %
holds:
\begin{align}\label{eq:Condition**}
\begin{aligned}
&\text{every $(x,z)$ satisfying $A^jx^j + B^jf_j(z)\in\K_j$ for all $j\in[p]$}\\
&\qquad\text{also satisfies $A^j\alpha_j x^j + B^jf_j(z)\in \K_j$ for any $\alpha_j\geq 1$ for all $j\in[p]$.}
\end{aligned}\tag{$\star\star$}
\end{align}
\end{proposition}
\proof{{\bf Proof.%
}}
	Similar to the proof of Proposition~\ref{Condition*}.
\Halmos
\endproof
In Section~\ref{sec:soc-multiple}, inspired by fractional binary programming, we study a generalization of $\cH$ and $\cR$ that takes \emph{multiple} conic quadratic constraints, whose convex hull can be characterized based on Corollary~\ref{cor:multipleFunctions}.

There are  several functions $f$ for which  $\conv(\Epi(f))$ is known explicitly; we list some of these examples below. 
\begin{itemize}
\item \emph{Submodular functions.} \cite{lovasz1983submodular} characterizes the description of $\conv(\Epi(f))$ when $f$ is a submodular function (see also \cite[Theorem 1]{atamturk2008polymatroids}). In fact, we observe that submodularity is the common structure that arises in the applications of our theoretical developments so far.
Due to this, we provide a detailed discussion on submodular functions and their epigraphs in Appendix~\ref{sec:inequalities}. In particular, we discuss the class of inequalities that describe the convex hull of the epigraph of a submodular function, known as the \emph{polymatroid} inequalities. %
(See Section \ref{sec:applications} for the application of these inequalities in our setting.) %
Additionally, we highlight that \emph{multiple} conic constraints that use \emph{common binary variables} can be handled \emph{jointly} when the underlying functions are all submodular, based on the convex hull results of \cite{baumann2013submodular,Kilinc-Karzan2019joint-sumod}. The same result is unlikely to hold in general.%

\item \emph{Univariate convex or separable convex functions.} \cite{miller2003tight} give the convex hull of the epigraph of any function that is univariate convex or piecewise linear convex. %
In particular, the authors show that for the case of piecewise linear convex functions with non-integer breakpoints, mixed-integer rounding inequalities are sufficient to describe the convex hull. These results are applicable in our setting if we have $f(z)=g(\sum_{i=1}^n z_i)$  or $f(z)=\sum_{i=1}^n g_i(z_i)$ for convex $g, g_i, i\in[n]$. (See Example \ref{example}  for the application of these inequalities in our setting.) 

\item \emph{Low dimensional problems.} For small dimensional problems, the convex hull description can be obtained by disjunctive cuts \citep{BALAS19793}.

\item \emph{Bisubmodular functions.} \cite{YU20215} provide the convex hull description of the epigraph of any bisubmodular function. Bisubmodular functions and their epigraphs are defined over $z\in\{0,\pm1\}^n$, not $z\in\{0,1\}^n$. Nevertheless,
the adaptation of the statement and proof of Theorem~\ref{thm:conv-trick} to cover the case when $\cG\subseteq\R^p_+\times\{0,\pm1\}^n$ is straightforward. 
\end{itemize}

We close this section by noting that the structure of the set $\cS(f,\K)$ allows us to easily embed constraints on the continuous variables $x$ as well.
\begin{remark}[Additional constraints on continuous variables]\label{rem:+conicConstraints}
Let $\C$ be a closed  convex pointed cone, %
and consider the set along with its transformation given by
\begin{align*}
&\left\{(x,z)\in\R^m\times\{0,1\}^n:~\exists y\in\R \text{ s.t. } y\geq f(z),~ \tilde{A}x+\tilde{B}y\in\tilde{\K},~C x\in\C \right\} \\
& =\left\{(x,z)\in\R^m\times\{0,1\}^n:~\exists y\in\R \text{ s.t. }y\geq f(z),~ Ax+ B y\in\K \right\},
\end{align*}
where we set $\K= \tilde{\K}\times\C$, $A=[\tilde{A}; C]$, and $B=[\tilde{B}; 0]$. \fullrank{(Also, note that if $\tilde A$ %
	has full column rank, so does 
	$A$.%
	) }Thus,
through this representation, we deduce that the additional conic constraints on the continuous variables $x$ can easily be embedded into our desirable form of the set $\cS(f,\K)$.
\ifx\flagJournal\true \epr \fi
\end{remark}

%% file: extension.tex
\subsection{Extensions}\label{sec:preliminaries}

Our results are applicable in the cases where the conic constraint is non-homogeneous, or the function $f(z)$ does not satisfy the nonnegativity assumption, or where we have only a partial convex hull description of the epigraph of $f(z)$  available.

\begin{remark}[Non-homogeneous conic constraints]\label{rem:nonhomogeneous}
Our results are still of interest when 
$\cS(f,\K)$ has a non-homogeneous constraint, i.e., $Ax +By+C\in \K$ %
for some $C\neq0$ instead of $Ax +By\in \K$. %
Indeed, by adding a new variable and an affine constraint, we can always rewrite the associated set $\cQ(\cG(f),\K)$  using a homogeneous conic constraint. That is,
\begin{align}\label{non-homogeneous}
\begin{aligned}
&\left\{(x,y,z)\in\R^m\times \R \times\{0,1\}^n:~~y\geq f(z),~ Ax +By+C\in \K   \right\}\\
&\qquad = \left\{(x,y,z)\in\R^m\times \R \times\{0,1\}^n:~~ \exists v\in\R \text{ s.t. } y\geq f(z),~ Ax +By+Cv\in \K, ~v=1   \right\}.
\end{aligned}
\end{align}
Here, $\left\{(x,y,v,z)\in\R^m\times\R\times\R\times\{0,1\}^n:y\geq f(z),~ Ax +By+Cv\in \K\right\}$ is of the form $\cQ(\cG(f),\K)$, and the set in~\eqref{non-homogeneous} is obtained from the intersection of this set and an affine hyperplane defined by $v=1$ after projecting out $v$. Therefore, %
we apply Theorem~\ref{thm:multipleFunctions} to obtain 
\begin{equation}\label{convex-relaxation}
\left\{(x,y,z)\in\R^m\times\R\times[0,1]^n:~\exists v\in\R\text{ s.t. } (y,z)\in\conv(\Epi(f)),~Ax +By+Cv\in \K, ~v=1\right\}
\end{equation} 
as a valid convex relaxation for the set in~\eqref{non-homogeneous}. 
\ifx\flagJournal\true \epr \fi
\end{remark}
While Remark~\ref{rem:nonhomogeneous} is useful, due to the presence of the affine constraint $v=1$ in~\eqref{convex-relaxation}, ~\eqref{convex-relaxation} may not provide a convex hull description of the set in~\eqref{non-homogeneous}. This is demonstrated in the following example.

\begin{example}\label{example}
In \eqref{non-homogeneous}, let $m=0$ (or $A=0$), $n=1, B=1, \K=\R_+, f(z)=\max\{0.5-z,\,0.5z-0.25\}$ and $C=-1/3$. Then, \eqref{non-homogeneous} is equivalent to
\begin{align*}
&\left\{(y,z)\in\R \times\{0,1\}:~y\geq 0.5-z,~ y\geq 0.5z-0.25,~ y-1/3\ge 0   \right\},
\end{align*}
and its convex hull is given by adding the  inequality $2y+z/3\ge 1$  to the continuous relaxation:%
\begin{align*}
\begin{aligned}
&\conv\left(\left\{(y,z)\in \R \times\{0,1\}:~~y\geq f(z),~ y-1/3\ge0   \right\}\right)\\
&\qquad = \left\{(y,z)\in\R \times[0,1]:~~ y\geq f(z),~ y-1/3\ge0,~ 2y+z/3\ge 1 \right\}.
\end{aligned}
\end{align*}

In contrast, if we follow the procedure described in Remark~\ref{rem:nonhomogeneous}, we will first obtain the convex hull of the epigraph of $f$, which is given by
\begin{align}\label{eq:exampleConv}
\begin{aligned}
&\!\!\!\!\!\!\!\!\conv\left(\left\{(y,z)\in \R \times\{0,1\}:y\geq f(z) \right\}\right)
 = \conv\left(\left\{(y,z)\in\R \times\{0,1\}: y\geq 0.5-z, ~y\geq 0.5z-0.25 \right\}\right)\\
&\quad = \left\{(y,z)\in\R \times[0,1]:~  y\geq 0.5-z, ~y\geq 0.5z-0.25,~2y+0.5z\ge 1 \right\}.
\end{aligned}
\end{align}
Based on this, the relaxation given in \eqref{convex-relaxation} becomes
\begin{align}\label{eq:exampleRelax}
&\left\{(y,z)\in\R \times[0,1]:~y\geq 0.5-z, ~y\geq 0.5z-0.25,~ y\ge 0,~ y-1/3\ge0,~2y+0.5z\ge 1   \right\}.
\end{align}
Clearly, the sets defined in \eqref{eq:exampleConv} and \eqref{eq:exampleRelax} are different; see Figure~\ref{fig:example}.
\begin{figure}
\begin{center}
\begin{tikzpicture}[scale = 2.5]%
\draw[ step = 0.25cm, black, dotted, very thin] ( -0.9 ,-0.6 ) grid ( 1.5, 1.3 );
\draw[black,<->, thick] ( 0, 1.3 ) node[ anchor = north east ] {\large{$y$}} -- ( 0, -0.6 ) ;
\draw[black,<->, thick] ( 1.5,0 ) node[ anchor = south east ] {\large{$z$}} -- ( -1.1, 0 ) ;
\foreach \y in {0,1,...,1} 
\draw[black] (-.1pt,\y) -- (-.2pt,\y) node[anchor = north east] {\small{\y}};
\foreach \x in {-1,1} 
\draw[black] (\x,.1pt) -- (\x,-.2pt) node[anchor = north west] {\small{\x}};

 \fill[black,%
 fill opacity=1] (0,0.5) circle (1pt);
 \fill[black,%
 fill opacity=1] (1,1/3) circle (1pt);

\draw[blue, thick] ( -0.5, 1/3 ) -- ( 1.5,1/3 ); 
\draw[blue, thick]  ( -0.7, 1/3 ) node[anchor = north] {$y\ge1/3$};

\draw[blue, thick] ( -0.5, 1 ) -- ( 1,-0.5 ); 
\draw[blue, thick]  ( -0.5, 1 ) node[anchor = south] {$y\ge0.5-z$};

\draw[blue, thick] ( -0.5, -0.5 ) -- ( 1.5,0.5 ); 
\draw[blue, thick]  ( -0.5, -0.5 ) node[anchor = east] {$y\ge0.5z-0.25$};

\draw[blue, very thick, dashed] ( -0.5, 3.5/6 ) -- ( 1.5,0.25 ); 
\draw[blue, thick]  ( -0.5, 3.5/6 ) node[anchor = east] {$2y+ z/3\ge1$};

\draw[red, very thick, dashed] ( -0.5, 1.25/2 ) -- ( 1.5,1/8 ); 
\draw[red, thick]  ( -0.5, 1.25/2 ) node[anchor = south east] {$2y+ 0.5z\ge1$};
				
\fill[blue!30!white, fill opacity=0.5] (0,1.3) -- (0,0.5) -- (1,1/3) -- (1, 1.3) -- cycle;
				
\end{tikzpicture}
\caption{Illustration of Example~\ref{example}.}\label{fig:example}
\end{center}
\end{figure}
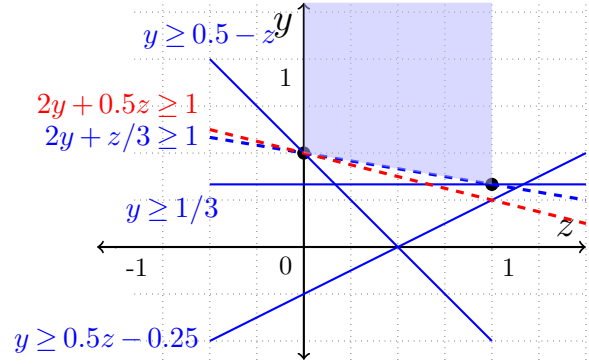
\ifx\flagJournal\true \epr \fi 
\end{example}

\begin{remark}[%
set functions]\label{rem:nonnegativity}
Consider the set $\overline\cS(h,\K)$ where $h$ is neither nonnegative nor nonpositive. Let $h_{\min}\coloneqq \min_{z\in\{0,1\}^n} h(z)$. 
%\revise{Here, $h_{\min}$ is finite as $\{0,1\}^n$ is finite.}  \sk{[Not sure, if finiteness of $z$ suffices, e.g., $h(z)=-1/(\sum_i z_i)$.]} \fkk{[I agree that we cannot directly argue that $h_{\min}$ is finite.]} Then, 
Assuming that $h_{\min}$ is finite,  the function $f(z)\coloneqq h(z)-h_{\min}$ is a nonnegative function and we have
\begin{align*}
\overline\cS(h,\K) &= \left\{(x,z)\!\in\R^m\!\times\!\{0,1\}^n\!\!:\!%
Ax +Bh(z)\in \K
 \right\}%
= \left\{(x,z)\!\in\R^m\!\times\!\{0,1\}^n\!\!:\!%
 Ax +Bf(z)+B h_{\min}\in \K
   \right\} \!.
\end{align*}
Now, we can apply the transformations from Remark~\ref{rem:nonhomogeneous} and use Corollary~\ref{cor:singleFunction} to arrive at the desired representation.
\ifx\flagJournal\true \epr \fi
\end{remark}

\begin{remark}[Supermodular functions]\label{rem:supermodularTransformation}
Suppose the function of interest is supermodular. Let $h:\{0,1\}^n\to\R$ be a supermodular function and $\K$ be a closed  convex pointed cone. %
If $h$ is nonnegative, then one would want to apply  Corollary~\ref{cor:singleFunction} on the set $
\overline\cS(h,\K) = \left\{(x,z)\in\R^m\times\{0,1\}^n:~ Ax +Bh(z)\in \K   \right\}$ and convexify it by introducing an auxiliary variable $y$ and replacing $Ax +Bh(z)\in \K$ with the constraints $(y,z)\in\conv(\epi(h))$ and $Ax+By\in\K$. %
However,  the description of %
$\conv(\epi(h))$ when $h$ is supermodular is not easy. Nevertheless, we can still utilize our technique of transforming the homogeneous conic constraint into a non-homogeneous one and applying Remark~\ref{rem:nonhomogeneous}. To this end, define $h_{\max}\coloneqq\max_{z\in\{0,1\}^n}h(z)$. Assuming $h_{\max}$ is finite, let $f(z)\coloneqq-h(z)+h_{\max}$. 
%\revise{Note that $h_{\max}$ is finite as $\{0,1\}^n$ is finite.} \sk{[Not sure, if finiteness of $z$ suffices, e.g., $h(z)=1/(\sum_i z_i)$.]}  \fkk{[I agree that we cannot directly argue that $h_{\max}$ is finite.]} 
Then, we know that $f$ is a nonnegative submodular function, and we have %
\begin{align*}
\overline\cS(h,\K) &= \left\{(x,z)\in\R^m\times\{0,1\}^n:~~ Ax +B(h_{\max}-f(z))\in \K   \right\}\\
&= \left\{(x,z)\in\R^m\times\{0,1\}^n:~~ \exists v\in\R \text{ s.t. } Ax +Bh_{\max}v-Bf(z)\in \K,~v=1 \right\}.
\end{align*}
Therefore, as $f(z)$ is submodular, we can add the extended polymatroid inequalities for $f$ to strengthen the set $\overline \cS(h,\K)$. 

We note that this transformation applied to a supermodular function is indeed easy as computing $h_{\max}$ amounts to minimizing a submodular function, which can be done in polynomial time.
\ifx\flagJournal\true \epr \fi	
\end{remark}

Whether the transformation for a supermodular function to a submodular function, discussed in Remark~\ref{rem:supermodularTransformation}, is useful or not depends on the structure of the cone $\K$ and the parameters $A,B$ used in the definition of the set $\overline\cS(h,\K)$. For example, consider any supermodular function $h(z)$ and the simple set defined by 
$\overline\cS(h,\R_+) \coloneqq \left\{(x,z)\in\R\times\{0,1\}^n:~ x - h(z)\in \R_+   \right\}$, which is nothing but the epigraph of the supermodular function $h(z)$ and perhaps %
its convex hull does not admit an easy characterization. If we apply the transformation described in Remark~\ref{rem:supermodularTransformation}, we end up with the set
$\overline\cS(h,\R_+) =\left\{(x,z)\in\R\times\{0,1\}^n:~ \exists v\in\R \text{ s.t. } x - h_{\max}v+f(z)\in \R_+, ~v=1   \right\}$ where $f$ is the submodular function defined by $\coloneqq-h(z)+h_{\max}$. While this set fits the form that we are interested in, after implementing this transformation, unfortunately we cannot apply the necessary reformulation step described in Proposition~\ref{Condition*} for this transformed set as any $(x,v,z)\in\R^2\times\{0,1\}^n$ satisfying $f(z)\geq h_{\max}v-x$ will not necessarily satisfy $f(z)\geq \alpha (h_{\max}v-x)$ for all $\alpha\ge1$. 
That said, there are cases when this remark is of interest. In particular, we illustrate in Section~\ref{sec:BSS} how the transformation from Remark~\ref{rem:supermodularTransformation} can be applied %
the context of Best Subset Selection. %

\begin{remark}\label{rem:polarInequalities}
For general (not necessarily submodular) functions $f$, \citet{atamturk2020polar} introduced the class of polar inequalities that are valid for $\widehat{\cS}(f,\K)$ for any nonnegative function $f$, which can be used to  strengthen the continuous relaxation; see Appendix~\ref{sec:inequalities}. 
\ifx\flagJournal\true \epr \fi
\end{remark} 

%% file: application.tex
\section{Applications}\label{sec:applications}

In this section, we present several optimization problems in which sets of the form $\cS(f,\K)$ appear as a substructure. In this respect, we highlight two forms of objective functions that immediately lead to our desired structure:
(i) $\min~ \sqrt{f(z)+ \|Dx+d\|_2^2}$,
(ii) $\min~ {\|Dx+d\|_2^2 \over f(z)}$,
where $f(z)$ is a nonnegative function, in addition to norm constraints of the form $\|(z;Dx)\|\le t$ which can be expressed as cone constraints involving closed convex pointed cones.

In Section~\ref{sec:soc}, we discuss the implications of our results in the context of applications from \cite{atamturk2019conic-quad} which have the objectives of form (i) where $f$ is submodular, and in Section~\ref{sec:soc-multiple} we explore their use in the case of fractional programming which involves several conic constraints on common binary variables (note that convex hull result for such cases were not provided in the literature previously). In Section~\ref{sec:BSS}, we point out connections with variants of best subset selection problem~\citep{gomez2020fractional} in which the problems have objectives of the form (ii) and the function $f$ is supermodular. %
Finally, in Section~\ref{sec:misdp}, we highlight how our framework can be exploited within the context of sparse approximation of the sum of positive semidefinite matrices.

\input{SOC_submodular}

\subsection{Best subset selection}\label{sec:BSS}

\cite{gomez2020fractional}  study the following general model for the \emph{Best Subset Selection} (BSS) problem:
\begin{align}\label{eq:BSS}
\min_{\beta\in\R^n, z\in\{0,1\}^n} \left\{ ~\frac{\|a-U\beta\|_2^2}{g(\sum_{i\in[n]}z_i)}:~ - M z_i \leq \beta_i \leq M z_i,~ \forall i\in[n] \right\}, 
\end{align}
where $a\in\R^k,~U\in\R^{k\times n}$ are data, $M\in\R_+$ corresponds to big-M value, and $g:\R_+\to\R_+$ is a non-increasing convex function.
The best subset selection problem in linear regression is to find a sparse subset of regressors that best fits the data. To this end, in Problem~\eqref{eq:BSS}, the regression variables are $\beta$ and each binary variable $z_i$ determines whether a regressor $\beta_i$ is selected. While mean squared error (MSE) is a popular criterion to measure the goodness of fit, other criteria such as Akaike Information Criterion (AIC), corrected AIC (AICc), and Bayesian  Information Criterion (BIC) have also been proposed. The latter three criteria have desirable properties that address some shortcomings of the MSE criterion, but they involve (non-convex) logarithmic terms in the objective  function that call for advanced solution methods. In particular, for AIC and BIC, we have $g(\sum_{i\in[n]}z_i)=e^{-\alpha \sum_{i\in[n]}z_i}$, and for AICc we have  $g(\sum_{i\in[n]}z_i)=e^{-2\alpha/(\alpha- \sum_{i\in[n]}z_i)}$ for an appropriate choice of $\alpha\ge 0$. We refer the reader to \cite{gomez2020fractional} for more details.

\cite{gomez2020fractional} work with the following reformulation of~\eqref{eq:BSS}:
\begin{align}\label{eq:BSS-re}
\min_{\beta\in\R^n, z\in\{0,1\}^n, s\in\R_+} \left\{ ~\frac{\|a-U\beta\|_2^2}{s}:~ s\leq g(\sum_{i\in[n]}z_i),~- M z_i \leq \beta_i \leq M z_i,~ \forall i\in[n] \right\}.
\end{align}
In~\eqref{eq:BSS-re}, $s\leq g(\sum_{i\in[n]}z_i)$ is equivalent to $-s\geq -g(\sum_{i\in[n]}z_i)$, and since $-g(\sum_{i\in[n]}z_i)$ is submodular, the feasible region of~\eqref{eq:BSS-re} can be strengthened by applying the extended polymatroid inequalities for $-g$. However, the objective of~\eqref{eq:BSS-re} is a still \emph{fractional function}. Consequently, \cite{gomez2020fractional}  apply a customized Newton-type method after \emph{parameterizing} the fraction.

In contrast, we observe here that~\eqref{eq:BSS} admits a SOC reformulation with a linear objective and whose feasible region contains a substructure that fits our framework. Define $h(z)\coloneqq g(\sum_{i\in[n]}z_i)$. By introducing a new variable $t\in\R_+$ to capture the objective function, we observe that \eqref{eq:BSS} is equivalent to the following problem:
\begin{align*}
\min_{t\in\R_+,\beta\in\R^n, z\in\{0,1\}^n} \left\{ t:~~t\cdot h(z) \geq \|a-U\beta\|_2^2, ~ - M z_i \leq \beta_i \leq M z_i,~ \forall i\in[n] \right\}. 
\end{align*}
Note that the nonlinear constraint $t\cdot h(z) \geq \|a-U\beta\|_2^2$ in this formulation is equivalent to requiring 
\[
[2a-2U\beta;\ t-h(z);\ t+h(z)]\in\L^{k+2}.   
\]
Furthermore, when we define the function $h(z)\coloneqq g(\sum_{i\in[n]}z_i)$ based on a nonnegative non-increasing convex function $g$, we deduce that $h(z)$ is nonnegative and supermodular. By using the transformation from Remark~\ref{rem:supermodularTransformation}, we define $h_{\max}\coloneqq\max_{z\in\{0,1\}^n}\{h(z)\}$  (note that since $h$ is a composition of a non-increasing function with an affine function, $h_{\max}=h(0)=g(0)\ge 0$, and for the MSE, AIC, BIC, and AICc criteria, $h_{\max}<\infty$) and $f(z)\coloneqq-h(z)+h_{\max}$, and arrive at the equivalent conic constraint
\[
[2a-2U\beta;\ t+f(z)-h_{\max};\ t-f(z)+h_{\max}]\in\L^{k+2},   
\]
where $f$ is a nonnegative submodular function.
In order to homogenize this constraint as we did in Remark~\ref{rem:nonhomogeneous}, we introduce another decision variable $v\in\R$ and the constraint $v=1$. By letting $x=[t;v;\beta]$ and $m=n+2$, our conic constraint becomes $\tilde{A}x + \tilde{B}f(z)\in\L^{k+2}$, where we select 
\[
\tilde{A}=\begin{bmatrix} 0 & 2a & -2U \\ 
1 & -h_{\max} & 0 \\
1 & ~~h_{\max} & 0
\end{bmatrix} 
\quad \text{and}\quad 
\tilde{B}=\begin{bmatrix}  0_k \\ 1 \\ -1 \end{bmatrix}.
\] 
Hence, we arrive at the equivalent problem 
\begin{align*}
\min_{x\in\R^m, z\in\{0,1\}^n} \left\{x_1:~\tilde{A}x + \tilde{B}f(z)\in\L^{k+2},~x_1\in\R_+,~x_2=1,~ - M z_i \leq x_{i+2} \leq M z_i,~ \forall i\in[n] \right\}, 
\end{align*}
where %
$f(z)=-g(\sum_{i\in[n]} z_i) +\max_{\bar{z}\in\{0,1\}^n} \{g(\sum_{i\in[n]} \bar{z}_i)\}$ is a nonnegative submodular function.
Consider a point $(x,z)$ that satisfies $\tilde Ax+\tilde B f(z) \in \Le^{k+2}$. Then, $x=[t;v;\beta]$ and  $t\cdot (h_{\max}v-f(z)) \geq \|a-U\beta\|_2^2$ holds. Note that this constraint along with $t\ge 0$ implies that $h_{\max}v-f(z)\ge0$. Moreover, since $f$ is a nonnegative function, we have $h_{\max}v-0\ge h_{\max}v-f(z)$. Then, $t\cdot (h_{\max}v-0) \geq t\cdot (h_{\max}v-f(z)) \geq \|a-U\beta\|_2^2$ holds, i.e., $x$ satisfies $\tilde Ax\in \Le^{k+2}$ as well. 
Then, Remark~\ref{rem:condition*} implies that Condition~\eqref{eq:Condition*} holds and we can thus apply Proposition~\ref{Condition*} to arrive at the equivalent formulation 
\begin{align}
\label{eq:SOCP:BSS}
\min_{x\in\R^m, y\in\R_+,z\in\{0,1\}^n} \left\{x_1:~ 
\begin{array}{l}
	y\geq f(z), ~\tilde{A}x + \tilde{B}y\in\L^{k+2},~x_1\in\R_+,~x_2=1 \\ [2ex]
	- M z_i \leq x_{i+2} \leq M z_i,~ \forall i\in[n]
\end{array} 
\right\}. 
\end{align}

Note that in \eqref{eq:SOCP:BSS} we still cannot drop the binary requirements  on $z$ due to the presence of big-M constraints. Thus, our results provide only a tighter convex relaxation. 
Our developments in this application highlight that one can exploit the submodularity structure in this problem all the while using the standard optimization solvers without the need to develop specialized algorithms such as Newton-type methods or parametrization of the fractional objective. In Appendix~\ref{app:numerical}, we provide a numerical study on this problem to illustrate this point.

\input{misdp.tex}

%% file: SOC_submodular.tex
\subsection{Recovering the results of \cite{atamturk2019conic-quad}}\label{sec:soc}

Recall the following sets studied by \cite{atamturk2019conic-quad}:
\begin{align*}
\cH&=\left\{(x,z)\in\R^{m}_+\times\{0,1\}^n:~\sqrt{\sigma+\sum_{i\in[n]}c_iz_i+ \sum_{j\in[m-1]} d_j x_j^2 }\leq x_{m} \right\},\\
\cR&=\left\{(x,z)\in\R^{m}_+\times\{0,1\}^n:~\sigma+\sum_{i\in[n]}c_iz_i+ \sum_{j\in[m-2]} d_j x_j^2 \leq 4x_{m-1}x_{m}   \right\} .
\end{align*}

Consider the function $f:\{0,1\}^n\to\R_+$ defined by $f(z)\coloneqq \sqrt{\sigma+\sum_{i\in[n]}c_iz_i}$. Note that $f$ is submodular, provided that $\sigma+\sum_{i\in[n]}c_iz_i$ is nonnegative for $z\in \{0,1\}^n$. Hence, by Theorem~\ref{thm:lovasz}, $\conv(\Epi(f))$ is described by the extended polymatroid inequalities for $f-\sqrt{\sigma}$ and $\mathbf{0}\leq z\leq \mathbf{1}$. Motivated by this observation, \citet{atamturk2019conic-quad}  prove the following results %
by analyzing KKT conditions of a generic linear optimization problem over the domain $\cH$ or $\cR$. 
\begin{proposition}[{\citet[Proposition~5]{atamturk2019conic-quad}}]\label{conic-quadratic}
Let $f(z)\!=\!\sqrt{\sigma+\sum_{i\in[n]}c_iz_i}$ where $\sigma+\sum_{i\in[n]}c_iz_i\geq 0$ for $z\in \{0,1\}^n$, and let $\cH$ be defined as in~\eqref{H}. %
Then
\begin{align*}
\conv(\cH)
&=\left\{(x,z)\in\R^m_+\times\R^n:~\exists y\in\R_+~~\textup{s.t.}~(y,z)\in\conv(\Epi(f)),~ \sqrt{y^2+ \sum_{i\in[m-1]} d_i x_i^2} \leq x_m\right\} .
\end{align*}
\end{proposition}
\begin{proposition}[{\citet[Proposition~6]{atamturk2019conic-quad}}]\label{rotated-conic-quadratic}
	Let $f(z)\!=\!\sqrt{\sigma+\sum_{i\in[n]}c_iz_i}$ where $\sigma+\sum_{i\in[n]}c_iz_i\geq 0$ for $z\in \{0,1\}^n$, and let $\cR$ be defined as in~\eqref{R}. Then %
	\begin{align*}
	\conv(\cR)
	&=\left\{(x,z)\in\R^m_+\times[0,1]^n:~\exists y\in\R_+~~\textup{s.t.}~(y,z)\in\conv(\Epi(f)),~ {y^2+ \sum_{i\in[m-2]} d_i x_i^2} \leq 4x_{m-1}x_m\right\}.
	\end{align*}
\end{proposition}

We next show that these results are simple corollaries of Theorem~\ref{thm:multipleFunctions}. Recall $\L^{k+1}=\{(\xi,y)\in\R^k\times \R:~y\geq \norm{\xi}_2\}$ denotes the SOC in $\R^{k+1}$ (note that $\L^{k+1}$ is a closed convex pointed cone).  Throughout, let $e_i$ be the unit vector of appropriate dimension with a one in the $i$th component, and let $\Diag(\cdot)$ represent a diagonal matrix with the specified diagonal entries. 

\begin{corollary}\label{cor:cone+nonnegativeFunc}
	Theorem~\ref{thm:multipleFunctions} implies Proposition~\ref{conic-quadratic}.
\end{corollary}
\proof{{\bf Proof.%
}}
Note that the mixed-integer set $\cH$ considered in Proposition~\ref{conic-quadratic} satisfies
\begin{align*}	
\cH %
&=\left\{(x,z)\in\R^{m}_+\times\{0,1\}^n:~\sqrt{(f(z))^2+ \sum_{j\in[m-1]} d_j x_j^2 }\leq x_{m} \right\} \\
&=\left\{(x,z)\in\R^{m}_+\times\{0,1\}^n:~\left[f(z);\ \sqrt{d_1} x_1;\ \ldots;\ \sqrt{d_{m-1}} x_{m-1};\ x_m\right]\in \Le^{m+1} \right\} \\
&=\left\{(x,z)\in\R^{m}_+\times\{0,1\}^n:~\tilde Ax+\tilde B f(z) \in \Le^{m+1} \right\},
\end{align*}
where $\tilde B\coloneqq e_1\in\R^{m+1}$ and $\tilde A\coloneqq [0^\top;\ \Diag(\sqrt{d_1};\ \ldots;\ \sqrt{d_{m-1}};\ 1)]\in\R^{(m+1)\times m}$. Note that if $(x,z)$ satisfies $\tilde Ax+\tilde B f(z) \in \Le^{m+1}$, then due to the structure of $\tilde A,\tilde B$ and the cone $\Le^{m+1}$, we deduce that $x$ satisfies $\tilde Ax\in \Le^{m+1}$ as well. 
Then, by Remark~\ref{rem:condition*}, Condition~\eqref{eq:Condition*} holds. Therefore, by Proposition~\ref{Condition*}, 
\[
\cH=\left\{(x,z)\in\R^{m}_+\times\{0,1\}^n:~\exists y\in\R \text{ s.t. } y\geq f(z),~\tilde Ax+\tilde B y \in \Le^{m+1} \right\}.
\]
\fullrank{Note also that the matrix $[\tilde{B},\tilde{A}]$ is diagonal and thus $\tilde A$ %
	has full column rank.}
The result then follows by applying 
Theorems~\ref{thm:multipleFunctions}--\ref{thm:edmonds}, Corollary~\ref{cor:multipleSubmodularFunctions},  and Remark~\ref{rem:+conicConstraints} so that we can handle the nonnegativity constraint on the continuous variables $x\in\R^m_+$ (which is nothing but a very simple conic constraint) in addition to the SOC constraint $\tilde Ax+\tilde B y \in \Le^{m+1}$.
\Halmos
\endproof

\begin{corollary}\label{cor:rotated}
	Theorem~\ref{thm:multipleFunctions} implies Proposition~\ref{rotated-conic-quadratic}.
\end{corollary}
\proof{{\bf Proof.%
}}
Note that $y^2+ \sum_{i\in[m-2]} d_i x_i^2 \leq 4x_{m-1}x_m$ is  a rotated SOC constraint given by
\[
\left[y; \sqrt{d_1}x_1; \ldots; \sqrt{d_{m-2}}x_{m-2};\, x_{m-1}-x_m;\, x_{m-1} + x_m \right]^\top \in \L^{m+2}.
\]
By defining  $\tilde B\coloneqq e_1\in\R^{m+2}$ and $\tilde A\coloneqq [0^\top;\ \Diag(\sqrt{d_1};\ \ldots;\ \sqrt{d_{m-2}});\ 0^\top;\ 0^\top ]+e_{m}[0,\ldots,0,1,-1] + e_{m+1}[0,\ldots,0,1,1]\in\R^{(m+2)\times m}$, %
this constraint is equivalent to $\tilde{A}x +\tilde{B}y \in\L^{m+2}$. 
Then,
\begin{align*}
\cR %
&=\left\{(x,z)\in\R^{m}_+\times\{0,1\}^n:~(f(z))^2+ \sum_{j\in[m-2]} d_j x_j^2 \leq 4x_{m-1}x_{m}   \right\} \\
&=\left\{(x,z)\in\R^{m}_+\times\{0,1\}^n:~\tilde{A}x +\tilde{B}f(z) \in\L^{m+2} \right\} \\
&=\left\{(x,z)\in\R^{m}_+\times\{0,1\}^n:~\exists y\in\R \text{ s.t. } y\geq f(z), ~\tilde{A}x +\tilde{B}y \in\L^{m+2}   \right\} ,
\end{align*}
where the last line follows from Proposition~\ref{Condition*}  (exactly as in the case of Corollary~\ref{cor:cone+nonnegativeFunc}). 
Furthermore, 
by applying the transformation from Remark~\ref{rem:+conicConstraints} so that we can handle the nonnegativity constraint on the continuous variables $x\in\R^m_+$ in addition to the SOC constraint $\tilde Ax+\tilde B y \in \Le^{m+2}$, 
we conclude that the result follows from Theorem~\ref{thm:multipleFunctions}, Corollary~\ref{cor:multipleSubmodularFunctions}\fullrank{(note again that the matrix $\tilde A$ %
	has full column rank)}, and Theorems~\ref{thm:lovasz}~and~\ref{thm:edmonds}. %
\Halmos
\endproof

\subsection{Multiple conic quadratic constraints}\label{sec:soc-multiple}
Recall that Theorem~\ref{thm:multipleFunctions} may take \emph{multiple} functions into account at the same time. Based on Theorem~\ref{thm:multipleFunctions}, we can characterize the convex hull of a set defined by multiple conic quadratic constraints, generalizing the results on $\cH$ and $\cR$. For two finite disjoint sets $H,R$ of indices, the following set is defined by $|H|$ conic quadratic constraints of the type used in $\cH$ and $|R|$ constraints of the type used in $\cR$:
\begin{equation}\label{multiple-conic-quadratic}
\cM:=\left\{
(x,z)
:~~
\begin{aligned}
&(x,z)\in\R^{m(|H|+|R|)}_+\times\{0,1\}^n,\\
&\sqrt{\sigma_\ell+\sum_{i\in[n]}c_{\ell,i}z_i+ \sum_{j\in[m-1]} d_{\ell,j} x_{\ell,j}^2 }\leq x_{\ell,m},~~\forall \ell\in H\\
&\sigma_\ell+\sum_{i\in[n]}c_{\ell,i}z_i+ \sum_{j\in[m-2]} d_{\ell,j} x_{\ell,j}^2 \leq 4x_{\ell,m-1}x_{\ell,m},~~\forall \ell\in R
\end{aligned} \right\}.
\end{equation}
We are interested in $\conv(\cM)$. As in the proofs of Corollary~\ref{cor:cone+nonnegativeFunc} and~\ref{cor:rotated}, we can rewrite $\cM$ as
\[
\cM:=\left\{(x,z)\in\R^{m(|H|+|R|)}_+\times\{0,1\}^n:~~
\begin{aligned}
&\tilde A^\ell x^\ell + \tilde B^\ell f_\ell(z)\in\L^{m+1},~~\forall \ell\in H\\
&\tilde A^\ell x^\ell + \tilde B^\ell f_\ell(z)\in\L^{m+1},~~\forall \ell\in R
\end{aligned} \right\},
\]
where $x^\ell = (x_{\ell,1},\ldots, x_{\ell,m})^\top$ and $\tilde A^\ell, \tilde B^\ell$ are defined as in the proofs of Corollary~\ref{cor:cone+nonnegativeFunc} and~\ref{cor:rotated} for $\ell\in H\cup R$. Notice that $\cM$ is of the form $\overline{\cS}(\{f_j\}_{j\in[p]},\{\K_j\}_{j\in[p]})$ defined as in~\eqref{eq:Sbar-multiple}. Moreover, as both $\cH$ and $\cR$ satisfy Condition~\eqref{eq:Condition*}, $\cM$ satisfies the condition~\eqref{eq:Condition**}. Therefore, we can apply Corollary~\ref{cor:multipleFunctions} and Proposition~\ref{Condition**} to obtain the following proposition characterizing the convex hull of $\cM$.

\begin{proposition}\label{pr:multiple-conic-quadratic}
For $\ell\in H\cup R$, let $f_\ell(z)=\sqrt{\sigma_\ell+\sum_{i\in[n]}c_{\ell,i}z_i}$ where $\sigma_\ell+\sum_{i\in[n]}c_{\ell,i}z_i\geq0$ for $z\in\{0,1\}^n$, and let $\cM$ be defined as in~\eqref{multiple-conic-quadratic}. Then
\[
\conv(\cM)=\left\{
(x,z)
:~~
\begin{aligned}
&(x,z)\in\R^{m(|H|+|R|)}_+\times[0,1]^n,\\
&\exists y_\ell\in\R_+\text{ s.t. }(y_\ell,z)\in\conv(\epi(f_\ell)),~\sqrt{y_\ell^2+ \sum_{j\in[m-1]} d_{\ell,j} x_{\ell,j}^2 }\leq x_{\ell,m},~~\forall \ell\in H\\
&\exists y_\ell\in\R_+\text{ s.t. }(y_\ell,z)\in\conv(\epi(f_\ell)),~y_\ell^2+ \sum_{j\in[m-2]} d_{\ell,j} x_{\ell,j}^2 \leq 4x_{\ell,m-1}x_{\ell,m},~~\forall \ell\in R
\end{aligned} \right\}.
\]
\end{proposition}
We remark that Proposition~\ref{pr:multiple-conic-quadratic} generalizes Propositions~\ref{conic-quadratic} and~\ref{rotated-conic-quadratic}. For the rest of this section, we list some applications of Proposition~\ref{pr:multiple-conic-quadratic}.

\begin{remark}[fractional binary programs]
We can use $\cM$ to model optimization problems of the following form:
\begin{equation}\label{fractional}
\min_{z\in \cX}\sum_{\ell\in R}\frac{a_{\ell,0}+\sum_{i\in[n]}a_{\ell,i}z_i}{b_{\ell,0}+\sum_{i\in[n]}b_{\ell,i}z_i},
\end{equation}
where $\cX\subseteq\{0,1\}^n$ and $a_{\ell,0},a_{\ell,i},b_{\ell,0},b_{\ell,i}$ for $i\in[n]$ are all nonnegative numbers. The fractional optimization model~\eqref{fractional} is used in a wide range of application domains including modeling multinomial logit (MNL) choice models in assortment optimization, set covering, market share based facility location, stochastic service systems, bi-clustering, and optimization of boolean query for databases ~\cite[see,][and references therein]{fractional2017}. We can reformulate~\eqref{fractional} by introducing an auxiliary variable for each fraction. Note that
\[
\frac{a_{\ell,0}+\sum_{i\in[n]}a_{\ell,i}z_i}{b_{\ell,0}+\sum_{i\in[n]}b_{\ell,i}z_i}\leq u_\ell\quad\iff\quad a_{\ell,0}+\sum_{i\in[n]}a_{\ell,i}z_i\leq u_\ell v_\ell,~~ v_\ell=b_{\ell,0}+\sum_{i\in[n]}b_{\ell,i}z_i.
\]
Then~\eqref{fractional} is equivalent to
\[
\min_{z\in \{0,1\}^n,\ u,v\in\R^\ell}\left\{\sum_{\ell\in R}u_\ell:~~
\begin{aligned}
&4a_{\ell,0}+\sum_{i\in[n]}4a_{\ell,i}z_i\leq 4u_\ell v_\ell,~~\forall \ell\in R,\\
&v_\ell=b_{\ell,0}+\sum_{i\in[n]}b_{\ell,i}z_i,~~\forall \ell\in R
\end{aligned}
\right\}.
\]
In particular, %
$4a_{\ell,0}+\sum_{i\in[n]}4a_{\ell,i}z_i\leq 4u_\ell v_\ell$ for $\ell\in R$ give rise to a set of the form~$\cM$.
\ifx\flagJournal\true \epr \fi	
\end{remark}

%% file: misdp.tex
\subsection{Sparse sums of positive semidefinite matrices}\label{sec:misdp}

A classic result on approximate Carath\'{e}odory theorem due to \cite{althofer1994sparse} states that given vectors $u_1, \dots, u_n \in [0, 1]^n$, and $\alpha \in \R_+^n$ with $\sum_{i \in [n]} \alpha_i = 1$, for any $\varepsilon \in (0, 1)$, there exists a sparse vector $\beta \in \R^n_+$ with $\sum_{i \in [n]} \beta_i = 1$ and cardinality of $\cO(\log(n) / \varepsilon^2)$ such that $\| \sum_{i \in [n]} \beta_i u_i - \sum_{i \in [n]} \alpha_i u_i  \|_\infty \leq \varepsilon$. 
 A matrix generalization of the approximate Carath\'{e}odory theorem is presented in \cite{silva2015sparse}. That is, given  symmetric positive semidefinite matrices $U_1, \dots, U_n \in \S_+^N$ and $\alpha \in \R_+^n$ with $\sum_{i \in [n]} \alpha_i = 1$ such that $\sum_{i\in[n]}\alpha_i U_i$ is positive definite, for any $\varepsilon \in (0, 1)$, there exists a sparse vector $\beta \in \R^n_+$ with $\sum_{i \in [n]} \beta_i = 1$ with cardinality of at most $\cO(n / \varepsilon^2)$ such that
\begin{align}
\label{eq:sparse:sum}
	(1 - \varepsilon) \Big( \sum_{i \in [n]} \alpha_i U_i \Big) \preceq \sum_{i \in [n]} \beta_i U_i \preceq (1 + \varepsilon) \Big( \sum_{i \in [n]} \alpha_i U_i \Big).
\end{align}
This result has a number of interesting applications, such as graph sparsification, sparse solutions to semidefinite programs, and sparsification of unique games; see \cite{silva2015sparse} and references therein. %We next 
In the following we 
show that this problem fits into our framework of mixed-binary conic sets.

A simple variable substitution $U_i \gets (\sum_{i \in [n]} \alpha_i U_i)^{-\frac{1}{2}} U_i (\sum_{i \in [n]} \alpha_i U_i)^{-\frac{1}{2}}$ %
leads to the following equivalent representations of~\eqref{eq:sparse:sum}
\begin{align*}
	(1 - \varepsilon) I_N \preceq \sum_{i \in [n]} \beta_i U_i \preceq (1 + \varepsilon) I_N ~ \Longleftrightarrow ~ \max \Big\{ \lambda_{\max}\Big( \sum_{i \in [n]} \beta_i U_i - I_N \Big),~ \lambda_{\max}\Big( I_N - \sum_{i \in [n]} \beta_i U_i \Big) \Big\} \leq \varepsilon,
\end{align*}
where $\lambda_{\max}(A)$ denotes the largest eigenvalue of the symmetric matrix $A$ and $I_N$ is the $N\times N$ identity matrix.  
Since the goal is to find a sparse $\beta$ satisfying these relations, we arrive at the optimization problem
\begin{align*}
	\begin{array}{c@{\quad}l}
		\min\limits_{\beta, z} & \displaystyle \max \Big\{ \lambda_{\max}\Big( \sum_{i \in [n]} \beta_i U_i - I_N \Big),~ \lambda_{\max}\Big( I_N - \sum_{i \in [n]} \beta_i U_i \Big) \Big\} + \gamma f(z) \\[3ex]
		\text{s.t.} & \displaystyle \beta \in \R^n_+, ~ z \in \{0,1\}^n, ~ \sum_{i \in [n]} \beta_i = 1, ~\beta_i (1 - z_i) = 0, ~\forall i \in [n],
	\end{array}
\end{align*}
where the constant $\gamma$ is positive, the indicator variable $z$ controls the support of $\beta$, and the function $f$ is submodular. The penalty term $\gamma f(z)$ is added to the objective to promote sparsity. For instance, setting $f(z) = \sum_{i \in [n]} z_i$ yields the penalty term $\gamma \| \beta \|_0$. In general, $f$ incorporates prior knowledge or structural constraints on $\beta$; see, for example, \cite{bach2010structured}.

Using the simple epigraph reformulation technique and incorporating binary variables to model the sparsity penalty term results in the following problem
\begin{align*}
	\min_{ \beta \in \R_+^n,\, t \in \R_+,\, z \in \{0,1\}^n} \left\{ t:~~ 
	\begin{array}{l}
	\lambda_{\max}\left( \sum_{i \in [n]} \beta_i U_i - I_N \right) + \gamma f(z) \leq t, \\[2ex]
	\lambda_{\max}\left( I_N - \sum_{i \in [n]} \beta_i U_i \right) + \gamma f(z) \leq t,\\[2ex]
	\sum_{i \in [n]} \beta_i = 1,~~0 \leq \beta_i \leq z_i, ~\forall i \in [n]
	\end{array}
	\right\}.
\end{align*}
The nonlinear constraints $\lambda_{\max}( \sum_{i \in [n]} \beta_i U_i - I_N ) + \gamma f(z) \leq t$ and $\lambda_{\max}( I_N - \sum_{i \in [n]} \beta_i U_i ) + \gamma f(z) \leq t$ in this formulation are equivalent to requiring 
\begin{align*}
	\begin{cases}
		\sum_{i \in [n]} \beta_i U_i - I_N \preceq (t - \gamma f(z)) I_N \\[2ex]
    	I_N - \sum_{i \in [n]} \beta_i U_i \preceq (t - \gamma f(z)) I_N 	
	\end{cases}
	\, \Longleftrightarrow \quad
	0\preceq t I_{2N} + \sum_{i \in [n]} \beta_i \begin{bmatrix} -U_i & 0 \\ 0 & U_i \end{bmatrix} - \gamma f(z) I_{2N} + \begin{bmatrix} I_N & 0 \\ 0 & -I_N \end{bmatrix} .
\end{align*}
In order to homogenize the linear matrix inequality constraint, we introduce another decision variable $v \in \R$ and the constraint $v = 1$ as we suggested in Remark~\ref{rem:nonhomogeneous}. By letting $x = [t;v;\beta]\in\R\times\R\times\R^n$, our linear matrix inequality becomes $\sum_{i \in [n+2]} A_i x_i + B f(z)  \succeq 0$, where 
	\begin{align*}
		A_1 = I_{2N}, ~ A_2 = \begin{bmatrix} I_N & 0 \\ 0 & -I_N \end{bmatrix}, ~ B = -\gamma I_{2N}, ~ A_{i+2} = \begin{bmatrix} -U_i & 0 \\ 0 & U_i \end{bmatrix}, \quad \forall i \in [n].
	\end{align*} 
Hence, we arrive at the equivalent problem
\begin{align*}
	\min_{x \in \R_+^{2+n},\; z \in \{0,1\}^n} \left\{ x_1: \sum_{i \in [n+2]} A_i x_i + B f(z) \in \S_+^{2N}, \, x_2 = 1, \, \sum_{i \in [n]} x_{i+2} = 1, ~ 0 \leq x_{i+2} \leq z_i, ~\forall i \in [n]  \right\},
\end{align*}
where $\S_+^{2N}$ denotes the cone of ${2N}\times{2N}$ symmetric positive semidefinite matrices, which is a closed convex pointed cone.
Given that $B = -\gamma I_{2N}$ is a negative definite matrix, we immediately deduce that any point $(x,z)$ that satisfies $\sum_{i \in [n+2]} A_i x_i + B f(z) \in \S_+^{2N}$ will immediately also satisfy $\sum_{i \in [n+2]} A_i x_i \in \S_+^{2N}$. Thus, Remark~\ref{rem:condition*} implies that Condition~\eqref{eq:Condition*} holds and we can thus apply Proposition~\ref{Condition*} to arrive at the equivalent formulation
\begin{align*}
	\min_{x \in \R_+^{2+n},\; z \in \{0,1\}^n,\; y \in \R} \left\{ x_1:~ 
	\begin{array}{l}
		\sum_{i \in [n+2]} A_i x_i + B y \in \S_+^{2N}, \, x_2 = 1, \, \sum_{i \in [n]} x_{i+2} = 1 \\[2ex]
		0 \leq x_{i+2} \leq z_i, ~\forall i \in [n], ~y \geq f(z)
	\end{array} \right\}.
\end{align*}
This mixed-integer semidefinite problem can then be further strengthened by exploiting the submodularity of $f(z)$ via cut generation. The current off-the-shelf semidefinite programming solvers use interior point methods, which are not well-suited to warm-starting. Warm starting is a critical component needed for the success of cut generation and branch-and-bound procedures in integer programming, which in the context of this application would require solving numerous semidefinite problems from scratch, a prohibitively expensive task for the current solvers. Our theoretical results for the semidefinite case may bear fruit in practice once the mixed-integer semidefinite programming solvers reach the maturity to benefit from delayed cut generation. %

%% file: appendix.tex
{
\begin{APPENDICES}

\input{inequalities.tex}

\section{Numerical study on best subset selection problem}\label{app:numerical}

In this section, we study the best subset selection problem~\eqref{eq:BSS} for AIC, BIC and MSE criterion. We compare the performance of the fractional programming approach proposed in~\citep{gomez2020fractional} against our mixed-binary second-order conic programming (MISOCP) formulation presented in~\eqref{eq:SOCP:BSS}. The fractional programming approach is based on solving a series of parametric mixed-binary quadratic programming (MIQP) problems (where the objective is convex quadratic and the constraints are linear) and updating the corresponding parameter with a Newton step. In particular, for any given parameter $t$, each MIQP in~\citep{gomez2020fractional} is of the form
\begin{align}
	\label{eq:QP:BSS}
	\widehat d(t) := \min_{\beta\in\R^n, z\in\{0,1\}^n, s\in\R_+} \left\{ ~\|a-U\beta\|_2^2 - t \cdot s:~ s\leq g(\sum_{i\in[n]}z_i),~- M z_i \leq \beta_i \leq M z_i,~ \forall i\in[n] \right\}.
\end{align}
Let $(\beta^\star, z^\star, s^\star)$ be an optimizer of~\eqref{eq:QP:BSS}. Then, in the next step, the parameter $t$ is updated according to the step rule $t \gets \| a - U \beta^\star \| / g(z^\star)$ and the resulting MIQP is solved again.
This procedure is repeated until the terminal condition $\widehat d(t) \geq 0$ is met. \citep[Proposition~4]{gomez2020fractional} proves that this iterative scheme converges after at most $n + 1$ iterations. 

In the sequel, we refer to our formulation~\eqref{eq:SOCP:BSS} as ``MISOCP,'' and the fractional programming method established in~\citep{gomez2020fractional} as ``MIQP + FP.'' 
Note that the state-of-the-art MIQP solvers are currently more advanced and faster than MISOCP solvers, thanks to the use of a variant of the simplex method that allows for warm-starting when solving node relaxation problems. In contrast, the node relaxation problems in MISOCPs rely on interior point methods, which are not well-suited to warm-starting. This means that generating cuts and carrying out branch-and-bound procedures require solving many SOCP problems from scratch, which is an expensive task. We thus analyze the iterative procedure of \citep{gomez2020fractional} in which~\eqref{eq:QP:BSS} is formulated as an equivalent formulation MISOCP. We refer to this method as ``MISOCP + FP.''  
In all mixed-binary optimization problems, we first solve the ordinary least squares problem, and then set the big-$M$ constant to two times the $\ell_\infty$ norm of the solution. 
All experiments are run on an Intel XEON CPU with with an AMD Opteron 4184 processor with 12 CPUs and 70GB of RAM. All mixed-binary quadratic and second-order cone programs are solved with Gurobi 9.5 using the Gurobi interface in Python. In order to ensure the reproducibility of our experiments, we make all source codes available at \url{https://github.com/sorooshafiee/best\_subset\_selection}.

We conduct experiments on both real-world and synthetic datasets. We use the same real-world datasets used in~\citep{gomez2020fractional}; these are accessible at the UCI Machine Learning Repository~\citep{UCI}  
%~\citeappendix[]{UCI} 
except for ``Diabetes.''   We first use the AIC criteria, where the function $g(z) = \exp(-2 \sum_{i \in [n]} z_i / n)$. For these datasets and AIC criteria, we report in Table~\ref{table:AIC:comparison} the size of each dataset and the following statistics: the execution time in seconds, total number of user and lazy cuts, total number of branch and bound nodes, and the optimality gap when a time limit of one hour is reached. The number of iterations used in the fractional programming methods are reported in parentheses under the gap columns in this table. The optimality gap for our formulation is simply the gap provided by the solver, and for fractional programming approaches we compute and report the optimality gap  based on Equation~(27) of \citep{gomez2020fractional}.   

The results in Table~\ref{table:AIC:comparison} suggest that in the case of small size datasets, in terms of solution time, our MISOCP formulation is slower than the MIQP + FP method but it is significantly faster than the MISOCP + FP method. This is indeed due to the fact that the state of the art QP solvers are significantly faster than the SOCP solvers. As a result, the relaxation problems in the branch and bound algorithm are solved much faster in the MIQP formulation than in the MISOCP formulations. Moreover, we observed that the primal heuristics developed for MIQPs are far more advanced than the ones developed for MISOCPs. 
On the other hand, in the case of large datasets, our MISOCP formulation outperforms both fractional programming methods by a wide margin. In particular, the optimality gaps achieved by our MISOCP for the Crime, Diabetes, and Insurance datasets are significantly smaller than those obtained by the fractional programming methods.

\begin{table}[!tb]
\centering
	\resizebox{0.7\textwidth}{!}{
	\begin{tabular}{c|c|rrrr}
		 & &  & Total \# of & Total \# of &     \\ 
		Dataset                         & Method               & Time  & Cuts & Nodes & Gap\%    \\ \noalign{\vskip 1pt} \hline \noalign{\vskip 1pt} 
		\multirow{3}{*}{\shortstack{AutoMPG \\ $k=392, n=25$}}        & MISOCP           & 37.7  & 1805             & 3126              & 0.0        \\
		& MIQP + FP         & 18.8  & 9916             & 13133             & 0.0(4)      \\
		& MISOCP + FP       & 466.3 & 20552            & 38500             & 0.0(4)      \\ \noalign{\vskip 1pt} \hline \noalign{\vskip 1pt}
		\multirow{3}{*}{\shortstack{BreastCancer \\ $k=196, n=37$}}   & MISOCP           & 22.2  & 1595             & 2955              & 0.0        \\ 
		& MIQP + FP         & 3.1   & 1461             & 2115              & 0.0(2)      \\
		& MISOCP + FP       & 42.9  & 2951             & 5289              & 0.0(2)      \\ \noalign{\vskip 1pt} \hline \noalign{\vskip 1pt}
		\multirow{3}{*}{\shortstack{Crime \\ $k=1993, n=100$}}          & MISOCP           & 3600  & 10703            & 11942             & 14.2   \\
		& MIQP + FP         & 3600  & 120904           & 133040            & 107.1(2)    \\
		& MISOCP + FP       & 3600  & 7438             & 10498             & 571.4(1)    \\ \noalign{\vskip 1pt} \hline \noalign{\vskip 1pt}
		\multirow{3}{*}{\shortstack{Diabetes \\ $k=442, n=64$}}      & MISOCP           & 3600  & 150920           & 189293            & 10.2       \\
		& MIQP + FP         & 3600  & 680371           & 738238            & 21.2(2)     \\
		& MISOCP + FP       & 3600  & 85074            & 123185            & 72.3(2)     \\ \noalign{\vskip 1pt} \hline \noalign{\vskip 1pt}
		\multirow{3}{*}{\shortstack{Diabetes (Simple) \\ $k=442, n=10$}} & MISOCP           & 1     & 49               & 64                & 0.0      \\
		& MIQP + FP         & 0.3   & 70               & 67                & 0.0(2)      \\
		& MISOCP + FP       & 6.9   & 195              & 263               & 0.0(4)      \\ \noalign{\vskip 1pt} \hline \noalign{\vskip 1pt}
		\multirow{3}{*}{\shortstack{Flares \\ $k=1066, n=26$}}         & MISOCP           & 10.7  & 109              & 190               & 0.0          \\
		& MIQP + FP         & 0.9   & 53               & 39                & 0.0(2)      \\
		& MISOCP + FP       & 88.3  & 81               & 143               & 0.0(3)      \\ \noalign{\vskip 1pt} \hline \noalign{\vskip 1pt}
		\multirow{3}{*}{\shortstack{FlaresC \\ $k=442, n=10$}}        & MISOCP           & 18.7  & 232              & 441               & 0.0            \\
		& MIQP + FP         & 1.9   & 469              & 701               & 0.0(3)      \\
		& MISOCP + FP       & 113.9 & 555              & 1089              & 0.0(3)      \\ \noalign{\vskip 1pt} \hline \noalign{\vskip 1pt}
		\multirow{3}{*}{\shortstack{FlaresM \\ $k=442, n=10$}}        & MISOCP           & 3.9   & 41               & 75                & 0.0            \\
		& MIQP + FP         & 0.4   & 46               & 59                & 0.0(1)      \\
		& MISOCP + FP       & 15.2  & 39               & 87                & 0.0(1)      \\ \noalign{\vskip 1pt} \hline \noalign{\vskip 1pt}
		\multirow{3}{*}{\shortstack{FlaresX \\ $k=442, n=10$}}        & MISOCP           & 4     & 43               & 75                & 0.0           \\
		& MIQP + FP         & 0.6   & 51               & 43                & 0.0(2)      \\
		& MISOCP + FP       & 49    & 64               & 236               & 0.0(2)      \\ \noalign{\vskip 1pt} \hline \noalign{\vskip 1pt}
		\multirow{3}{*}{\shortstack{Housing \\ $k=506, n=13$}}        & MISOCP           & 8.5  & 116              & 190               & 0.0              \\
		& MIQP + FP         & 0.7   & 105              & 118               & 0.0(3)      \\
		& MISOCP + FP       & 12.2  & 325              & 515               & 0.0(3)      \\ \noalign{\vskip 1pt} \hline \noalign{\vskip 1pt}
		\multirow{3}{*}{\shortstack{Insurance \\ $k=5822, n=151$}}      & MISOCP           & 3600  & 9950             & 10937             & 5.0           \\
		& MIQP + FP         & 3600  & 24726            & 26600             & 55(1)       \\
		& MISOCP + FP       & 3600  & 3319             & 3626              & 82.4(1)     \\ \noalign{\vskip 1pt} \hline \noalign{\vskip 1pt}
		\multirow{3}{*}{\shortstack{Servo \\ $k=167, n=19$}}          & MISOCP           & 3.6   & 355              & 675               & 0.0             \\
		& MIQP + FP         & 0.7   & 664              & 861               & 0.0(3)      \\
		& MISOCP + FP       & 4.8   & 531              & 999               & 0.0(3)   
	\end{tabular}} \\
	\caption{AIC Performance of different mixed-binary optimization methods for real datasets. }
	\label{table:AIC:comparison}
\end{table}

We repeat the same experiment for the BIC criteria, where the function $g(z) = \exp(-\log(n) \sum_{i \in [n]} z_i / n)$. We summarize our result in Table~\ref{table:BIC:comparison}. In general, solving mixed-binary optimization problems with the BIC criteria appears to be more difficult than solving them with the AIC criteria. Nonetheless, we arrive at the same conclusion as in the previous case. Namely, in the case of small size datasets, in terms of solution time, our MISOCP formulation is slower to the MIQP + FP method but it is significantly faster than the MISOCP + FP method. On the other hand, in the case of large datasets, our MISOCP formulation outperforms both fractional programming methods by a wide margin. 

\begin{table}[!tb]
\centering
	\resizebox{0.7\textwidth}{!}{
	\begin{tabular}{c|c|rrrr}
		 & &  & Total \# of & Total \# of &     \\ 
		Dataset                         & Method               & Time  & Cuts & Nodes & Gap\%    \\ \noalign{\vskip 1pt} \hline \noalign{\vskip 1pt} 
		\multirow{3}{*}{\shortstack{AutoMPG \\ $k=392, n=25$}}        & MISOCP           & 81.7   & 3699             & 7207              & 0.0          \\
		& MIQP + FP         & 57.7   & 28436            & 37312             & 0.0(5)	    \\
		& MISOCP + FP       & 925.8  & 43992            & 85606             & 0.0(4)   	    \\ \noalign{\vskip 1pt} \hline \noalign{\vskip 1pt}
		\multirow{3}{*}{\shortstack{BreastCancer \\ $k=196, n=37$}}   & MISOCP           & 9.6    & 686              & 1320              & 0.0          \\ 
		& MIQP + FP         & 2.3    & 942              & 1403              & 0.0(2)    \\
		& MISOCP + FP       & 16.3   & 1090             & 2105              & 0.0(2)    \\ \noalign{\vskip 1pt} \hline \noalign{\vskip 1pt}
		\multirow{3}{*}{\shortstack{Crime \\ $k=1993, n=100$}}          & MISOCP           & 3600   & 18319            & 20565             & 16.5      \\
		& MIQP + FP         & 3600   & 136781           & 145403            & 1644.2(2) \\
		& MISOCP + FP       & 3600   & 7600             & 9413              & 2054.6(2) \\ \noalign{\vskip 1pt} \hline \noalign{\vskip 1pt}
		\multirow{3}{*}{\shortstack{Diabetes \\ $k=442, n=64$}}      & MISOCP           & 2949.9 & 86123            & 171532            & 0.0         \\
		& MIQP + FP         & 464    & 74383            & 93435             & 0.0(2)    \\
		& MISOCP + FP       & 3282.1 & 60174            & 118121            & 0.0(3)    \\ \noalign{\vskip 1pt} \hline \noalign{\vskip 1pt}
		\multirow{3}{*}{\shortstack{Diabetes (Simple) \\ $k=442, n=10$}} & MISOCP           & 0.9    & 41               & 50                & 0.0          \\
		& MIQP + FP         & 0.3    & 79               & 68                & 0.0(2)    \\
		& MISOCP + FP       & 6.8    & 179              & 247               & 0.0(4)    \\ \noalign{\vskip 1pt} \hline \noalign{\vskip 1pt}
		\multirow{3}{*}{\shortstack{Flares \\ $k=1066, n=26$}}         & MISOCP           & 7.2    & 83               & 109               & 0.0          \\
		& MIQP + FP         & 0.8    & 47               & 33                & 0.0(2)    \\
		& MISOCP + FP       & 188.4  & 226              & 452               & 0.0(4)    \\ \noalign{\vskip 1pt} \hline \noalign{\vskip 1pt}
		\multirow{3}{*}{\shortstack{FlaresC \\ $k=442, n=10$}}        & MISOCP           & 9.4    & 121              & 225               & 0.0         \\
		& MIQP + FP         & 1.1    & 219              & 325               & 0.0(2)    \\
		& MISOCP + FP       & 64.6   & 375              & 693               & 0.0(3)    \\ \noalign{\vskip 1pt} \hline \noalign{\vskip 1pt}
		\multirow{3}{*}{\shortstack{FlaresM \\ $k=442, n=10$}}        & MISOCP           & 3      & 21               & 35                & 0.0          \\
		& MIQP + FP         & 0.3    & 21               & 26                & 0.0(1)    \\
		& MISOCP + FP       & 7.1    & 20               & 35                & 0.0(1)    \\ \noalign{\vskip 1pt} \hline \noalign{\vskip 1pt}
		\multirow{3}{*}{\shortstack{FlaresX \\ $k=442, n=10$}}        & MISOCP           & 2.8    & 20               & 31                & 0.0          \\
		& MIQP + FP         & 0.3    & 25               & 27                & 0.0(1)    \\
		& MISOCP + FP       & 13.4   & 20               & 73                & 0.0(1)    \\ \noalign{\vskip 1pt} \hline \noalign{\vskip 1pt}
		\multirow{3}{*}{\shortstack{Housing \\ $k=506, n=13$}}        & MISOCP           & 14.7		   & 118              & 184               & 0.0          \\
		& MIQP + FP         & 0.7    & 147              & 162               & 0.0(3)    \\
		& MISOCP + FP       & 27.8   & 877              & 1341               & 0.0(5)   \\ \noalign{\vskip 1pt} \hline \noalign{\vskip 1pt}
		\multirow{3}{*}{\shortstack{Insurance \\ $k=5822, n=151$}}      & MISOCP           & 3600   & 7566             & 8519              & 4.2       \\
		& MIQP + FP         & 2573   & 7617             & 10221             & 0.0(2)    \\
		& MISOCP + FP       & 3600   & 3151             & 3439              & 1222.8(1) \\ \noalign{\vskip 1pt} \hline \noalign{\vskip 1pt}
		\multirow{3}{*}{\shortstack{Servo \\ $k=167, n=19$}}          & MISOCP           & 2.2    & 200              & 372               & 0.0          \\
		& MIQP + FP         & 0.5    & 430              & 545               & 0.0(3)    \\
		& MISOCP + FP       & 3.7    & 419              & 765               & 0.0(3)    
	\end{tabular}} \\
	\caption{BIC Performance of different mixed-binary optimization methods for real datasets. }
	\label{table:BIC:comparison}
\end{table}

We next consider the MSE criteria, where the function $g(z) = n - \sum_{i \in [n]} z_i$. In this case, the mixed-binary optimization problems cannot be strengthened using extended polymatroid inequalities due to the linearity of the function $g$. We summarize our result in Table~\ref{table:MSE:comparison}. In general, solving mixed-binary optimization problems with the MSE criteria is easier than solving them with the AIC and BIC criterion. Furthermore, we again observe that for large datasets, our MISOCP formulation outperforms both fractional programming methods by a wide margin. 

\begin{table}[!tb]
\centering
	\resizebox{0.7\textwidth}{!}{
	\begin{tabular}{c|c|rrrr}
		 & &  & Total \# of & Total \# of &     \\ 
		Dataset                         & Method               & Time  & Cuts & Nodes & Gap\%    \\ \noalign{\vskip 1pt} \hline \noalign{\vskip 1pt} 
		\multirow{3}{*}{\shortstack{AutoMPG \\ $k=392, n=25$}}        & MISOCP           & 1     & 0               & 39                & 0.0          \\
		& MIQP + FP         & 14.2  & 0             & 10295             & 0.0(5)    \\
		& MISOCP + FP       & 583.8 & 0            & 54476             & 0.0(6)    \\ \noalign{\vskip 1pt} \hline \noalign{\vskip 1pt}
		\multirow{3}{*}{\shortstack{BreastCancer \\ $k=196, n=37$}}   & MISOCP           & 120.3 & 0             & 16044             & 0.0          \\ 
		& MIQP + FP         & 8.9   & 0             & 6500              & 0.0(2)    \\
		& MISOCP + FP       & 243.5 & 0            & 28698             & 0.0(3)    \\ \noalign{\vskip 1pt} \hline \noalign{\vskip 1pt}
		\multirow{3}{*}{\shortstack{Crime \\ $k=1993, n=100$}}          & MISOCP         & 3600  & 0             & 4838              & 12.0         \\
		& MIQP + FP         & 3600  & 0           & 173831            & 1180.5(2) \\
		& MISOCP + FP       & 3600  & 0             & 10224             & 7960.8(1) \\ \noalign{\vskip 1pt} \hline \noalign{\vskip 1pt}
		\multirow{3}{*}{\shortstack{Diabetes \\ $k=442, n=64$}}      & MISOCP            & 3600  & 0           & 118454            & 17.9      \\
		& MIQP + FP         & 3600  & 0           & 811736            & 137.5(3)  \\
		& MISOCP + FP       & 3600  & 0            & 154539            & 2968.5(1) \\ \noalign{\vskip 1pt} \hline \noalign{\vskip 1pt}
		\multirow{3}{*}{\shortstack{Diabetes (Simple) \\ $k=442, n=10$}} & MISOCP        & 3.7   & 0               & 46                & 0.0          \\
		& MIQP + FP         & 0.2   & 0               & 30                & 0.0(2)    \\
		& MISOCP + FP       & 3     & 0               & 97                & 0.0(2)    \\ \noalign{\vskip 1pt} \hline \noalign{\vskip 1pt}
		\multirow{3}{*}{\shortstack{Flares \\ $k=1066, n=26$}}         & MISOCP          & 7.3   & 0                & 1                 & 0.0          \\
		& MIQP + FP         & 0.7   & 0               & 55                & 0.0(2)    \\
		& MISOCP + FP       & 119.1 & 0              & 366               & 0.0(2)    \\ \noalign{\vskip 1pt} \hline \noalign{\vskip 1pt}
		\multirow{3}{*}{\shortstack{FlaresC \\ $k=442, n=10$}}        & MISOCP           & 144.5 & 0              & 543               & 0.0          \\
		& MIQP + FP         & 0.8   & 0              & 354               & 0.0(2)    \\
		& MISOCP + FP       & 86    & 0              & 1119              & 0.0(2)    \\ \noalign{\vskip 1pt} \hline \noalign{\vskip 1pt}
		\multirow{3}{*}{\shortstack{FlaresM \\ $k=442, n=10$}}        & MISOCP           & 14.2  & 0              & 231               & 0.0          \\
		& MIQP + FP         & 0.7   & 0              & 112               & 0.0(2)    \\
		& MISOCP + FP       & 55.7  & 0              & 403               & 0.0(4)    \\ \noalign{\vskip 1pt} \hline \noalign{\vskip 1pt}
		\multirow{3}{*}{\shortstack{FlaresX \\ $k=442, n=10$}}        & MISOCP           & 15.3  & 0               & 87                & 0.0          \\
		& MIQP + FP         & 0.5   & 0               & 32                & 0.0(2)    \\
		& MISOCP + FP       & 27.3  & 0              & 212               & 0.0(2)    \\ \noalign{\vskip 1pt} \hline \noalign{\vskip 1pt}
		\multirow{3}{*}{\shortstack{Housing \\ $k=506, n=13$}}        & MISOCP           & 14    & 0              & 289               & 0.0          \\
		& MIQP + FP         & 0.6   & 0              & 171               & 0.0(3)    \\
		& MISOCP + FP       & 16.7  & 0              & 695               & 0.0(4)    \\ \noalign{\vskip 1pt} \hline \noalign{\vskip 1pt}
		\multirow{3}{*}{\shortstack{Insurance \\ $k=5822, n=151$}}      & MISOCP         & 3600  & 0              & 1056              & 5.9       \\
		& MIQP + FP         & 3600  & 0            & 27647             & 1409.7(1) \\
		& MISOCP + FP       & 3600  & 0             & 4379              & 1718(1)   \\ \noalign{\vskip 1pt} \hline \noalign{\vskip 1pt}
		\multirow{3}{*}{\shortstack{Servo \\ $k=167, n=19$}}          & MISOCP           & 4     & 0              & 637               & 0.0          \\
		& MIQP + FP         & 1.1   & 0             & 1383              & 0.0(4)    \\
		& MISOCP + FP       & 22.6  & 0             & 4008              & 0.0(7) 
	\end{tabular}} \\
	\caption{MSE Performance of different mixed-binary optimization methods for real datasets. }
	\label{table:MSE:comparison}
\end{table}

We now examine synthetic instances generated following the instance generation procedure from~\citep{bertsimas2016best}. For the dimensions $k=100$ and $n=1000$, we generate the design matrix $U$ from a multivariate Gaussian distribution with zero mean and covariance matrix $\Sigma$, where $\Sigma_{ij} = \rho^{|i-j|}$ for the autocorrelation parameter $\rho = 0.35$. The coefficients of the true regression parameter $\beta_0$ is set to be one for the first $b = \| \beta_0 \|_0$ elements and is zero otherwise. Finally, the response variable $y$ is generated from a multivariate Gaussian distribution with zero mean and covariance matrix $\sigma^2 I$, where we set the noise parameter $\sigma^2 := \beta_0^\top \Sigma \beta_0 / 2$ for different values of $\beta_0$. Note that this choice of $\sigma^2$ achieves the target signal to noise ratio level of $2$. 
In our experiments, we vary the cardinality of $\beta_0$ from 10 to 90, and generate 10 random instances at each level. We run the three approaches under a time limit of ten minutes and compare their performances in terms of average solution time, the optimality gap, total \# of cuts, and total \# of nodes in the B\&B tree. 
Figure~\ref{fig:AIC:comparison} summarizes our findings for the AIC criteria. For each approach (presented in a different color), the corresponding line provides the average performance metric for 10 instances and the shaded area corresponds to the range of the performance metric over these instances. We observe that when the $\| \beta_0 \|_0$ is small, the resulting mixed-binary optimization problems are easy, and the fractional programming methods converge to an optimal solution in a few number of iterations. As a result, the MIQP + FP method outperforms our MISOCP method. However, as we increase the cardinality of $\beta_0$, the optimization problems become more challenging, and the fractional programming methods require more iterations to find an optimal solution. Moreover, at each iteration of the fractional programming methods, the subproblems become increasingly more difficult to solve. Consequently, neither fractional programming method yields meaningful results when the time limit is reached. In contrast, our MISOCP formulation provides high quality solutions in this situation.
Additionally, we observe that the number of added cuts in the MIQP + FP method is much larger than the MISOCP methods. This is because Gurobi identifies incumbent solutions for MIQPs more quickly and easily than for MISOCPs. The same reasoning also applies to the number of nodes explored by the branch and bound algorithm.

We conclude this section by repeating the same experiment for the BIC and MSE criterion. Figures~\ref{fig:BIC:comparison} and \ref{fig:MSE:comparison} highlight our findings for the BIC and MSE criteria, which are consistent with our prior observation when using the BIC criterion. As no cut is added, the graph depicting the number of cuts for the MSE method is eliminated from Figure~\ref{fig:MSE:comparison}.

\begin{figure}[!tb]
    \center
	\includegraphics[width=0.41\columnwidth]{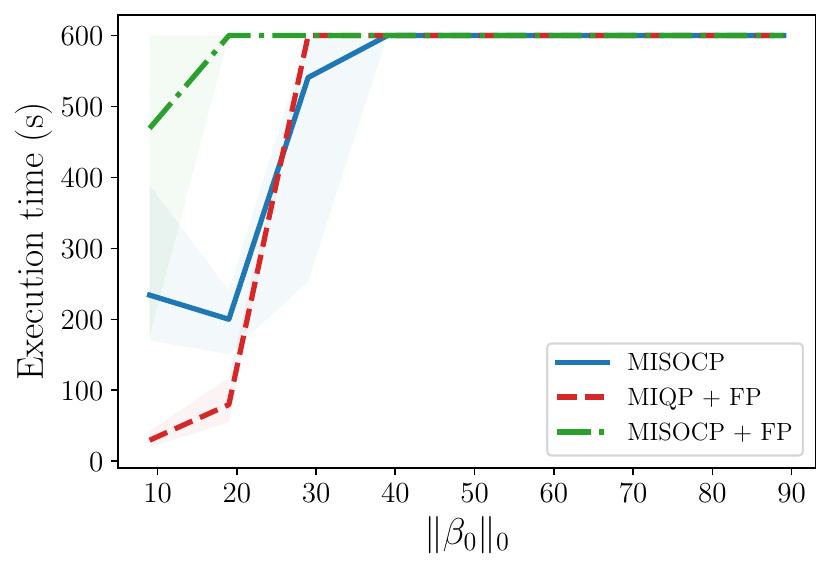}
	~\includegraphics[width=0.41\columnwidth]{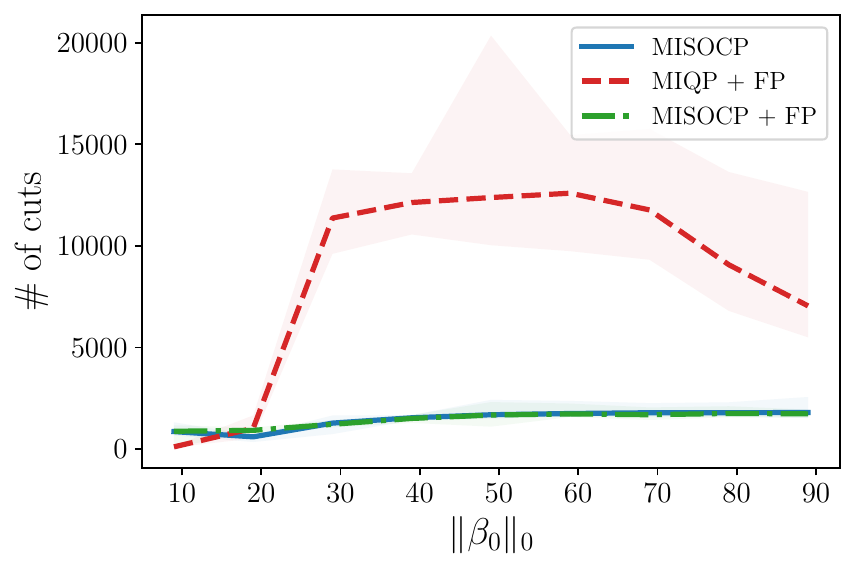} \\
	\includegraphics[width=0.41\columnwidth]{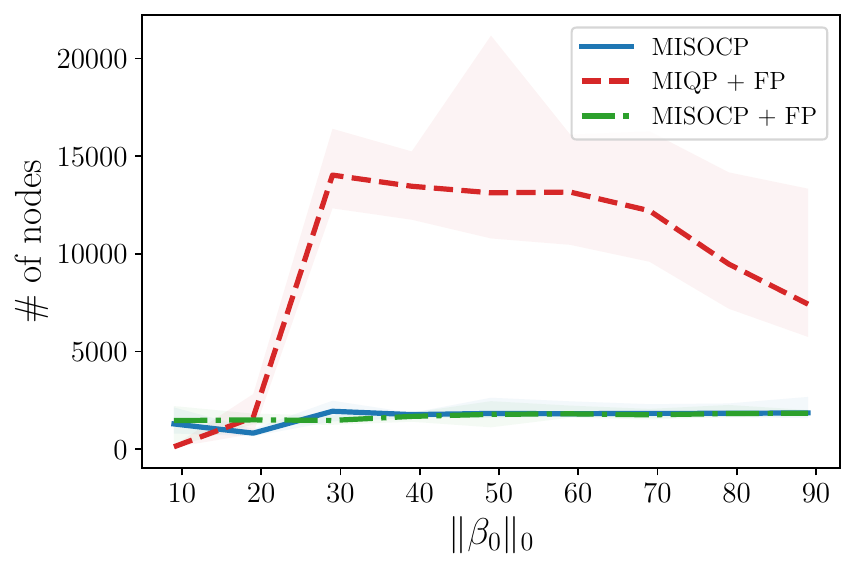}~
	\includegraphics[width=0.41\columnwidth]{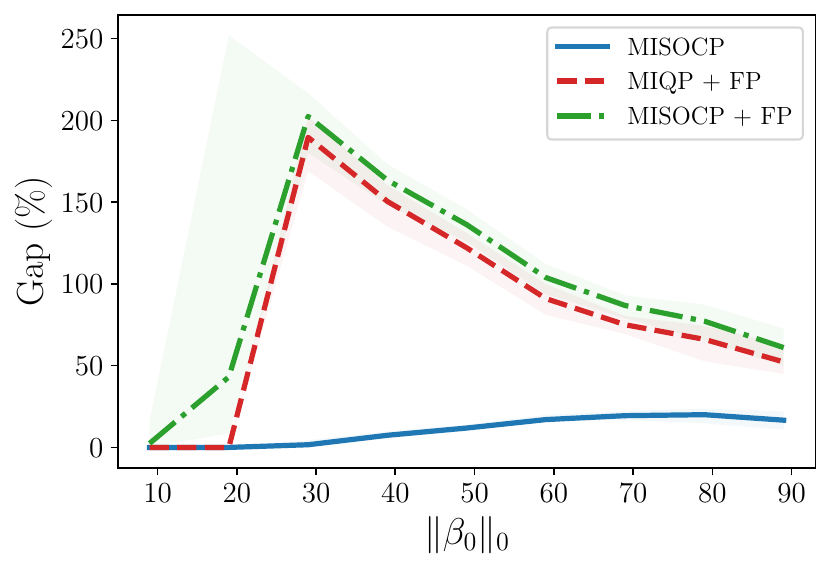}
	\caption{AIC Performance of different mixed-binary optimization methods for synthetic datasets. Solid lines (shaded regions) represent averages (ranges) across $10$ independent repetitions.}
	\label{fig:AIC:comparison}
\end{figure}

\begin{figure}[!tb]
    \center
	\includegraphics[width=0.41\columnwidth]{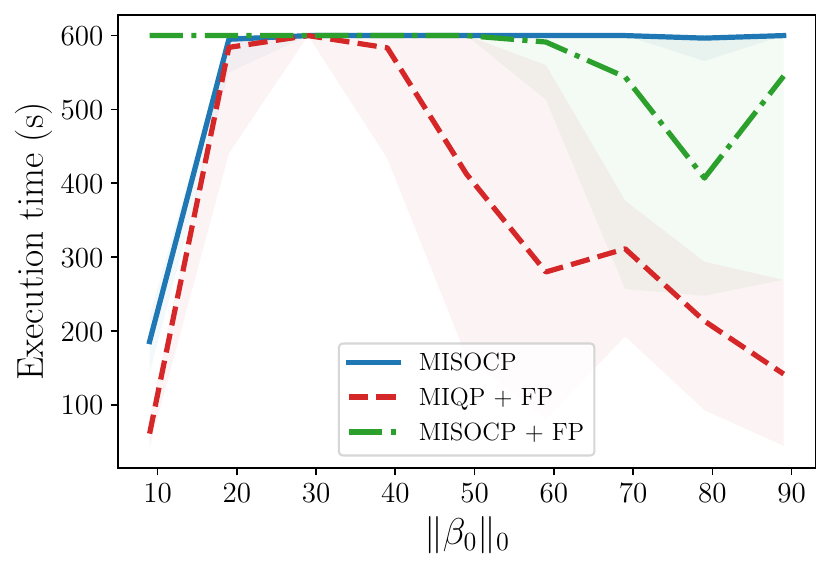}
	~\includegraphics[width=0.41\columnwidth]{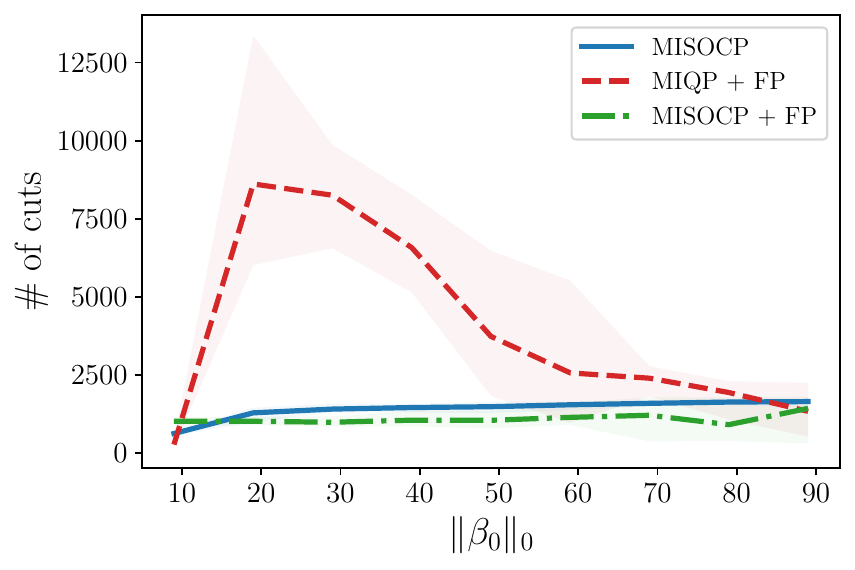} \\
	\includegraphics[width=0.41\columnwidth]{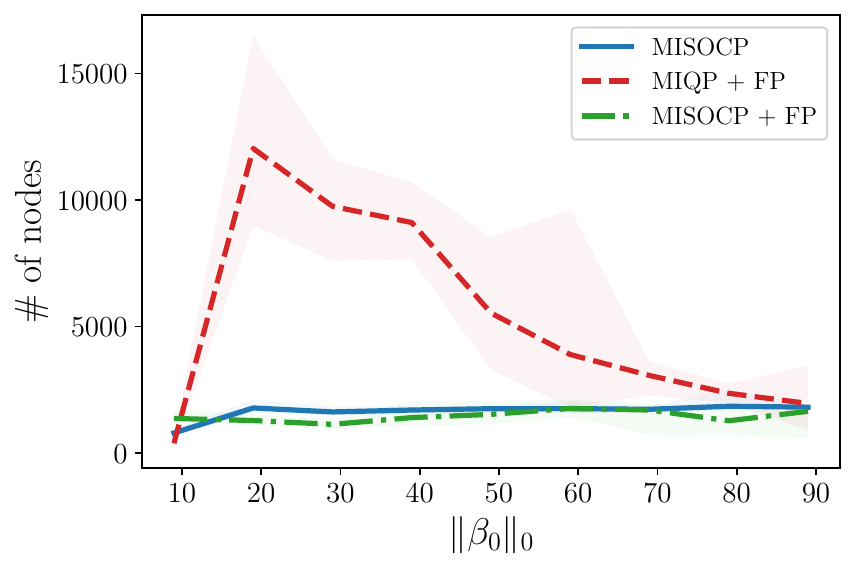}~
	\includegraphics[width=0.41\columnwidth]{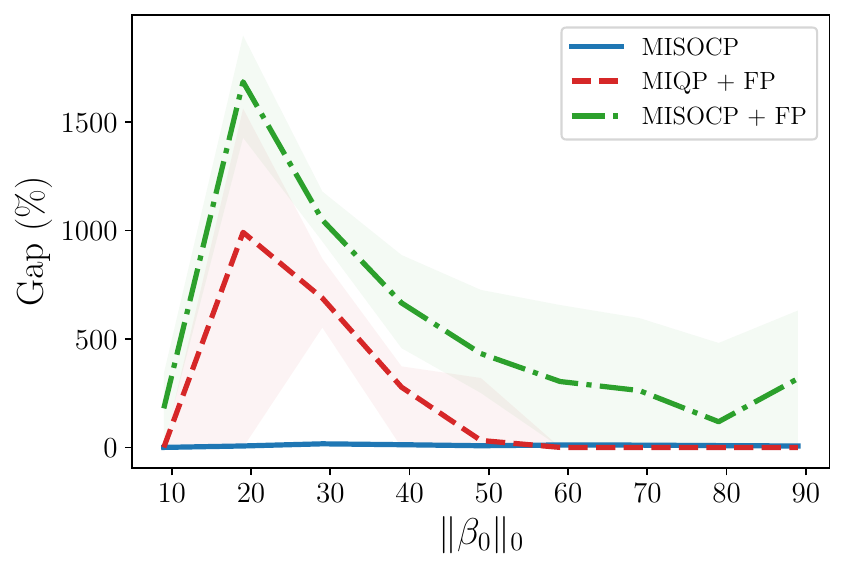}
	\caption{BIC Performance of different mixed-binary optimization methods for synthetic datasets. Solid lines (shaded regions) represent averages (ranges) across $10$ independent repetitions.}
	\label{fig:BIC:comparison}
\end{figure}

\begin{figure}[!tb]
    \center
	\includegraphics[width=0.41\columnwidth]{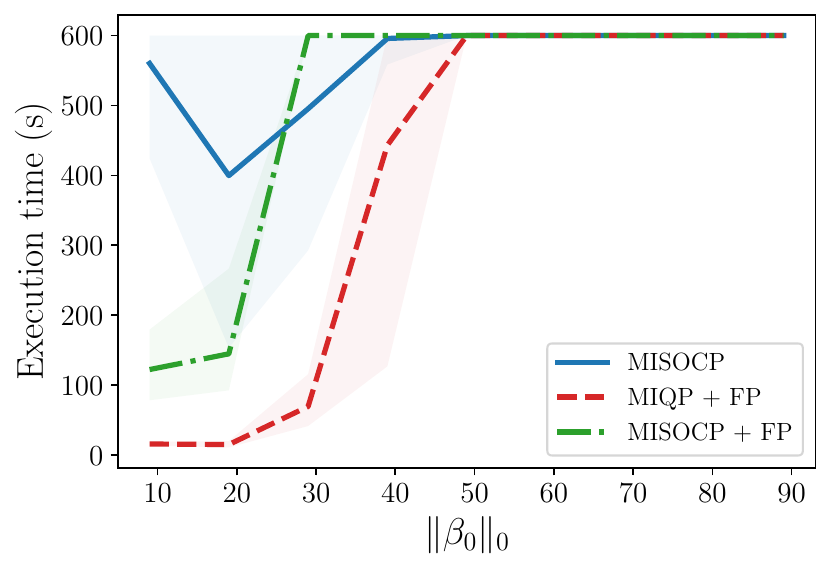}
	\\
	\includegraphics[width=0.41\columnwidth]{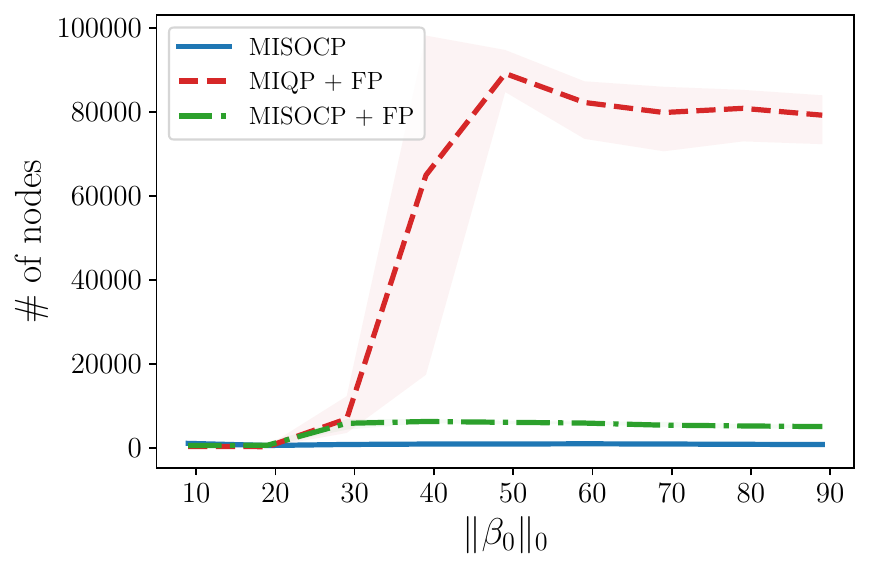}~
	\includegraphics[width=0.41\columnwidth]{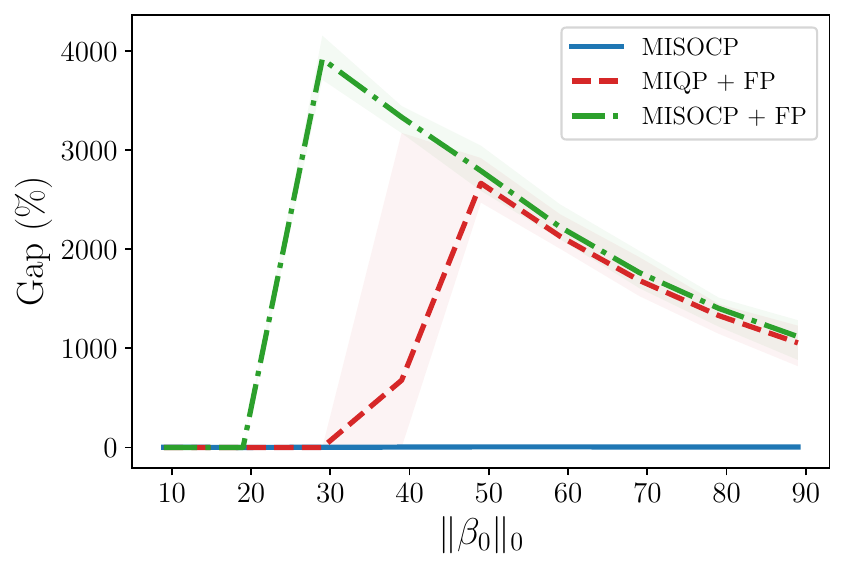}
	\caption{MSE Performance of different mixed-binary optimization methods for synthetic datasets. Solid lines (shaded regions) represent averages (ranges) across $10$ independent repetitions.}
	\label{fig:MSE:comparison}
\end{figure}

These results on real and synthetic instances collectively show that despite being a generic approach for conic optimization problems with binary variables and relying on standard optimization software, our technique is competitive (especially in large scale instances) against the state-of-the-art approach based on a Newton search procedure repeatedly solving optimization problems specifically designed for the Best Subset Selection problem. Notably, this is in fact the case despite the fact that the specialized algorithm exploits faster MIQP solvers as opposed to our approach relying on slower MISOCP solvers.

\input{droccp.tex}

\end{APPENDICES}
}

%% file: inequalities.tex
\section{Submodular functions and the extended polymatroid inequalities%
}\label{sec:inequalities}

In this appendix, %
we discuss two classes of inequalities, namely,  the extended polymatroid inequalities and the polar inequalities, which can be used to describe $\conv(\Epi(f))$ completely or partially (and thus  these give either $\conv(\cS(f,\K))$ or remain valid for it) when $f$ has desirable structure.

Given a set function $f:2^{[n]}\to\R$, where $2^{[n]}$ is the power set of $[n]$, let the \emph{associated polyhedron} of $f$ be defined as
\begin{equation}\label{ext-polymat}
\cP_f\coloneqq \left\{\pi\in\R^n:~\pi(V)\leq f(V),~\forall V\subseteq[n]\right\},
\end{equation}
where $\pi(V)\coloneqq \sum_{i\in V} \pi_i$ and $\pi(\emptyset)=0$. By slight abuse of notation, throughout, we  refer to  a set function  $f:2^{[n]}\to\R$  also as $f:\{0,1\}^n\to\R$, where 
$f(V)\coloneqq f(\mathbf{1}_V)$   for $V\subseteq[n]$ and $\mathbf{1}_V$ denotes the characteristic vector of $V$. When $f$ is a submodular function, i.e., $f$ satisfies
\[
f(S_1)+f(S_2)\geq f(S_1\cup S_2)+f(S_1\cap S_2),\quad \forall S_1,S_2\subseteq [n],
\]
$\cP_f$ is called the {\it extended polymatroid of $f$}. Note that $\cP_f$ is nonempty if and only if $f(\emptyset)\geq0$. In general, $f$ does not need to satisfy $f(\emptyset)\geq 0$. Nevertheless, we can take $f-f(\emptyset)$ instead so that $(f-f(\emptyset))(\emptyset)=0$. Hence, $\cP_{f-f(\emptyset)}$ is always nonempty. Hereinafter, we use notation $\tilde f$ to denote $f-f(\emptyset)$ for any set function $f$. In particular, $\tilde f$ is submodular when $f$ is submodular.

The associated polyhedron $\cP_f$ is instrumental in generating valid inequalities for $\epi(f)$ due to a close polarity relation between $\cP_f$ and $\conv(\epi(f))$. (We refer the reader to~\citep[Chapter I.4.5]{NWbook}  for a review of how the concept of polarity is used to obtain facets of a polyhedron.)  From this relation, \citep{atamturk2020polar} show that the so-called   \emph{polar inequalities}
\begin{equation}\label{polar-ineqs}
	y-f(\emptyset)\geq \pi^\top z,~\forall \pi\in \cP_{\tilde f}
	\end{equation}
are valid for $(y,z)\in \conv(\epi(f))$. 
\citep{atamturk2020polar}  also prove  that  inequalities \eqref{polar-ineqs} are facet-defining for $\conv(\epi(f))$ if and only if $\pi$ is an extreme point of $\cP_{\tilde f}$. Therefore,  the extreme points of $\cP_{\tilde f}$ characterize  facet-defining polar inequalities. In the case of a general set function $f$, the inclusion relationship in~\eqref{polar-ineqs} is strict  \citep[see Example 1,][]{atamturk2020polar}, which means that the polar inequalities may not be sufficient to describe the convex hull of $\epi(f)$.

When $f$ is submodular, the polar inequalities are precisely the \emph{extended polymatroid inequalities} \citep[]{atamturk2008polymatroids}. Moreover, it is well-known that when $f$ is submodular, the extended polymatroid inequalities indeed provide a complete description of $\conv(\epi(f))$. 
\begin{theorem}[{\citep{lovasz1983submodular}, \citep[Theorem 1]{atamturk2008polymatroids}}]\label{thm:lovasz}
	Let $f:\{0,1\}^n\to\R$ be a submodular function. Then, 
	\[
	\conv(\epi(f))=\left\{(y,z)\in\R\times[0,1]^n:~y-f(\emptyset)\geq \pi^\top z,~\forall \pi\in \cP_{\tilde f}\right\}.
	\]
\end{theorem}
\citet{edmonds1970polymatroid} provides the following explicit characterization of the extreme points of $\cP_{\tilde f}$.%
\begin{theorem}[\citep{edmonds1970polymatroid}]\label{thm:edmonds}
	Let $f:\{0,1\}^n\to\R$ be a submodular function. Then $\pi\in\R^n$ is an extreme point of $\cP_{\tilde f}$ %
	if and only if there exists a permutation $\sigma$ of~$[n]$ such that $\pi_{\sigma(t)}=f(V_{t})-f(V_{t-1})$, where $V_{t}=\{\sigma(1),\ldots,\sigma(t)\}$ for $t\in[n]$ and $V_{0}\coloneqq \emptyset$ by definition.
\end{theorem}
We note that when $f$ is not submodular, there are extreme points $\pi$ of $\cP_f$ are not necessarily of the form given in Therorem~\ref{thm:edmonds}. The proof of Theroem~\ref{thm:edmonds} yields an $O(n\log n)$ time algorithm for separating a violated extended polymatroid inequality given a point $(\bar y,\bar z)\in\R\times \R^n$, which amounts to solving $\max_{\pi}\left\{\bar z^\top \pi :~\pi\in \cP_{\tilde f}\right\}$~\citep[Section 2]{atamturk2008polymatroids}.

These results on submodular functions %
lead us to the following corollary of Theorem~\ref{thm:multipleFunctions}. 
\begin{corollary}\label{cor:generalization-submodular}
	Suppose $f:\{0,1\}^n\to\R_+$ is a nonnegative submodular function, $\K$ is a closed  convex pointed cone. %
	 \fullrank{and the matrix $A$ %
	 	has full column rank.}
	Then, $\conv(\cS(f,\K))$ is given by the extended polymatroid inequalities for $f-f(\emptyset)$ and the homogeneous conic constraint $Ax +By \in \K$. 
\end{corollary}
\proof{{\bf Proof. %
}}
	This is a direct consequence of Theorems~\ref{thm:multipleFunctions}~and~\ref{thm:lovasz}.
\Halmos
\endproof

\begin{corollary}\label{cor:multipleSubmodularFunctions}
For each $j\in[p]$, let $f_j:\{0,1\}^n\to\R_+$ be a nonnegative submodular function, $\K_j$ be a closed  convex pointed  cone. %
\fullrank{and the matrix $A^j$ %
	has full column rank for each $j\in[p]$.}
Then, the convex hull of $\cS(\{f_j\}_{j\in[p]},\{\K_j\}_{j\in[p]})$ defined in \eqref{eq:S-multiple} is described by $\widehat{\cS}(\{f_j\}_{j\in[p]},\{\K_j\}_{j\in[p]})$ as defined in \eqref{eq:Shat-multiple}, and moreover
\[
\widehat{\cS}(\{f_j\}_{j\in[p]},\{\K_j\}_{j\in[p]}) 
= \left\{(x,z)\in\R^{mp}\times[0,1]^{n}:~\exists y\in\R^p \text{ s.t. }~\begin{array}{l}(y_j,z)\in\conv(\Epi(f_j)),\\ A^jx^j + B^j y_j \in \K_j,\end{array}~\forall j\in[p]  \right\}. 
\]
\end{corollary}
\proof{{\bf Proof. %
}} 
Let 
$ %
\cG \coloneqq  \left\{(y,z)\in\R^{p}\times\{0,1\}^{n}:~y_j\geq f_j(z),~~\forall j\in[p]
\right\}.
$ %
Then, $\cG\subseteq\R^p_+\times\{0,1\}^n$ as $f_j$'s are nonnegative functions.  
The result follows from Theorem~\ref{thm:multipleFunctions} and the following fact. 
When $f_j$ is a nonnegative submodular function for each $j\in[p]$, from Theorem 2 of \citep{baumann2013submodular} \citep[see also, Proposition 1,][]{Kilinc-Karzan2019joint-sumod} we deduce that
\[
\conv (\cG) = \left\{(y,z)\in\R^{p}\times[0,1]^{n}:~(y_j,z)\in\conv(\Epi(f_j)),~~\forall j\in[p]
\right\}.
\Halmos
\]
\endproof

We refer the reader to \citep{edmonds1970polymatroid} and \citep{lovasz1983submodular} for a list of basic submodular functions as well as operations preserving submodularity.

\begin{remark}\label{rem:basicSubmodularF}
We close this section by emphasizing that our result holds not only for a function of the form $\sqrt{\sigma+\sum_{i\in[n]}c_iz_i}$ for $z\in\{0,1\}^n$, but also for general submodular functions, such as $g(\sigma+\sum_{i\in[n]}c_iz_i)$ where $g$ is concave. In particular, the constant elasticity of substitution function, $\left(\sum_{i\in[n]}c_i^pz_i^p\right)^{1/p}$ for any $c\in\R^n_+$ and $p\geq 1$ is submodular, and this property is used in Section~\ref{sec:DRO-CCP} when we consider  a norm constraint in the mixed-integer conic reformulation of a distributionally robust optimization problem. 
\ifx\flagJournal\true \epr \fi
\end{remark}

%% file: droccp.tex
\section{Distributionally robust chance-constrained programs under Wasserstein ambiguity} \label{sec:DRO-CCP}

Distributionally robust chance-constrained programming (DR-CCP) under Wasserstein ambiguity is formulated as 
\begin{equation}
\min\limits_{(x,z)}\left\{c^\top [x;z]:\ \sup_{\P \in \cF_N(\theta)} \P[\xi \not\in \cW(x,z)] \leq \epsilon,~ (x,z) \in \mathcal{X}\right\}.\label{eq:dr-ccp}
\end{equation}
Here, $c\in \R^{n+m}$ is a cost vector, $x\in \R^m$ is a vector of continuous decision variables, $z$ is a vector of $n$ binary decision variables, $\cX \subset \R^{n+m}$ is a compact domain for the decision variables, $\cW(x,z) \subseteq \R^K$ is a decision-dependent safety set, $\xi \in \R^K$ is a vector of  $K$ random variables with distribution $\P^*$, and $\epsilon \in (0,1)$ is the risk tolerance for the random variable $\xi$ falling outside the safety set $\cW(x,z)$.  Because the distribution $\P^*$  is often unavailable, independent and identically distributed (i.i.d.\@) samples $\{ \xi_i \}_{i \in [N]}$ are drawn from $\P^*$ to approximate  $\P^*$  using the empirical distribution $\P_N$ on these samples. To address the ambiguity in the true distribution, distributionally robust optimization model   \eqref{eq:dr-ccp} considers the worst-case probability of violating the safety constraints over a set of distributions on $\R^K$, given by   $\cF_N(\theta)$ , that contains the empirical distribution $\P_N$ where $\theta$ is a parameter that governs the size of the ambiguity set, and  the degree of conservatism of \eqref{eq:dr-ccp}.

\cite{chen2018data} and \cite{xie2018distributionally}   show that DR-CCP can be formulated as a mixed-integer conic program for certain $\cW(x,z)$, thus enabling their solution with standard optimization solvers. However, these MIP reformulations are  difficult to solve in certain cases for which the continuous relaxations provide weak lower bounds. \cite{lhs2020,rhs2020} consider a mixed-integer \emph{linear} substructure of the mixed-integer conic formulation of DR-CCP. The authors  propose  valid linear inequalities and other enhancements  that strengthen the continuous relaxation bounds and scale up the sizes of the problem that can be solved.

While previous research focused on the linear constraints, the strengthening of the mixed-integer conic reformulation of DR-CCP considered by \cite{xie2018distributionally} focuses on constraints of the form  
 $\|[\eta_1 x; \eta_1 z ; \eta_2]\|_*\le t$,  where $t$ is a continuous epigraph variable, $\eta_1,\eta_2\in\{0,1\}$ are constants with $\eta_1+\eta_2\ge 1$, and $\|\cdot\|_*$ is the dual of a norm $\|\cdot\|$. To obtain a strengthening from this conic constraint, we consider the case when $\eta_1=1$, i.e., when there is so-called left-hand side uncertainty. The norm is used to measure the Wasserstein distance and it is typical to consider an $\ell_q$-norm, whose dual norm is an $\ell_r$-norm with $\frac{1}{q}+\frac{1}{r}=1$. 
\cite{xie2018distributionally} pays specific attention to the pure-binary case, i.e.,  $m=0$ where the continuous variables $x$ are not present inside the norm. In this particular case,  \cite{xie2018distributionally} observes that the function $\|[\eta_1z ; \eta_2]\|_r$ is a submodular function and  the corresponding extended polymatroid inequalities can be applied to strengthen the set $\{(z,t):~\|[\eta_1z ; \eta_2]\|_r\leq t\}$ and this approach results in significant computational benefit. 
When the continuous variables $x$ are present, although $\|[\eta_1x; \eta_1 z ; \eta_2]\|_r$ is not even a set function. Nevertheless, our framework applies to constraints of the form $\|[\eta_1 x; \eta_1 z ; \eta_2]\|_r\le t$. 
To elaborate, let $\K_r^{m+2}$ denote the $r$-th order cone in $\R^{m+2}$ (note that $\K_r^{m+2}$ is a closed convex pointed cone), and define $f(z) \coloneqq \|[z ; \eta_2]\|_r$. Note that as $z_j^r=z_j$ for all $j\in[n]$ and $\eta_2^{1/r}=\eta_2$ as well, and thus from Remark~\ref{rem:basicSubmodularF} we deduce that  $f$ is a nonnegative submodular function. Moreover, 
  \begin{align*}
 &\left\{(x,t,z)\in\R^m_+\times\R_+\times\{0,1\}^n:~ \|[\eta_1x; \eta_1 z ; \eta_2]\|_r\le t\right\}\\
  & \quad =\left\{(x,t,z)\in\R^m_+\times\R_+\times\{0,1\}^n:~ \left[f(z);\  x_1;\ \ldots;\ x_{m};\ t\right]\in \K_r^{m+2} \right\} \\
  &\quad = \left\{(x,t,z)\in\R^m_+\times\R_+\times\{0,1\}^n:~\exists y\in\R \text{ s.t. }~ y\geq f(z), ~ \left[y;\  x_1;\ \ldots;\ x_{m};\ t\right]\in \K_r^{m+2} \right\} \\
  &\quad = \left\{(x,t,z)\in\R^m_+\times\R_+\times\{0,1\}^n:~\exists y\in\R \text{ s.t. }~y\geq f(z), ~\tilde A(x;t)+\tilde B y \in \K_r^{m+2} \right\},  
  \end{align*}
where  $\tilde B\coloneqq e_1\in\R^{m+2}$ and $\tilde A\coloneqq [0^\top;\ \Diag( 1;\ \ldots;\ 1)]\in\R^{(m+2)\times (m+1)}$. 
Here the second equation follows from Proposition~\ref{Condition*} because Condition~\eqref{eq:Condition*} holds 
(for any $(x,z)$ satisfying $\tilde Ax+\tilde B f(z) \in \K^{m+2}_r$, due to the structure of $\tilde A,\tilde B$ and the cone $\K^{m+2}_r$, we have $x$ satisfies $\tilde Ax\in \K^{m+2}_r$ as well, implying that we can apply Remark~\ref{rem:condition*} to deduce that Condition~\eqref{eq:Condition*} holds). 
\fullrank{Moreover, the matrix $[\tilde{B},\tilde{A}]$ is diagonal and thus $\tilde A$ %
	has full column rank, and} We can apply Corollary~\ref{cor:singleFunction} to obtain the convex hull of this set exactly.
Consequently, our results indicate that it is possible to exploit submodularity in the DR-CCP context even when  we have both continuous and binary decision variables. In particular, if $(x,z)$ is a mixed-binary decision vector, then by using Theorem~\ref{thm:multipleFunctions}, Corollaries~\ref{cor:singleFunction}~and~\ref{cor:multipleSubmodularFunctions}, and Theorems~\ref{thm:lovasz}~and~\ref{thm:edmonds} %
we can strengthen the resulting reformulation of DR-CCP under Wasserstein ambiguity.